\documentclass[12pt]{article}

\usepackage{indentfirst}
\usepackage{amsfonts}
\usepackage{amsmath}
\usepackage{amssymb}
\usepackage{amsbsy, amsthm}

\newtheorem{dfn}{Definition} [section]
\newtheorem{theorem}[dfn]{Theorem}
\newtheorem{lemma}[dfn]{Lemma}

\newtheorem{corollary}[dfn]{Corollary}
\newenvironment{pf}{\noindent{\bf Proof.}}
{\enspace\vrule height5pt depth0pt width5pt}

\def\X {{\mathcal X}}

\def\T {{\mathcal T}}

\def\Se {{\mathcal S}}

\def\mf {{\rm mf}}
\def\ins {{\rm ins}}
\def\C {{\mathcal C}}

\begin{document}

\centerline{{\large\bf EXCLUDING SUBDIVISIONS OF}}
\smallskip
\centerline{{\large\bf BOUNDED DEGREE GRAPHS}%
\footnote{Partially supported by NSF under Grant No.~DMS-1202640. 15 July 2014}%
\footnote{Revised: April 22, 2018. This material is based upon work supported by the National Science Foundation under Grant No.~DMS-1664593.}
}
\bigskip
\bigskip
\bigskip
\centerline{{\bf Chun-Hung Liu}%
\footnote{Current address: Department of Mathematics, Princeton University, Princeton, New Jersey 08544-1000, USA. 
Partially supported by NSF under Grant No.~DMS-1664593.}}
\smallskip
\centerline{and}
\smallskip
\centerline{{\bf Robin Thomas}%
}
\medskip
\centerline{School of Mathematics}
\centerline{Georgia Institute of Technology}
\centerline{Atlanta, Georgia  30332-0160, USA}
\bigskip

\begin{abstract}
\noindent
Let $H$ be a fixed graph. What can be said about graphs $G$ that have no subgraph
isomorphic to a subdivision of $H$?
Grohe and Marx proved that such graphs $G$ satisfy a certain structure theorem that
is not satisfied by graphs that contain a subdivision of a (larger) graph $H_1$.
Dvo\v{r}\'ak found a clever strengthening---his structure is not satisfied
 by graphs that contain a subdivision of a graph $H_2$, where $H_2$ has
``similar embedding properties" as $H$.
Building upon  Dvo\v{r}\'ak's theorem, we prove that said graphs $G$ satisfy
a similar structure theorem.
Our structure is not satisfied  by graphs that contain a subdivision of a graph $H_3$
that has similar embedding properties as $H$ and 
has the same maximum degree as $H$.
This will be important in a forthcoming application to well-quasi-ordering.
\end{abstract}

\section{Introduction}

In this paper {\em graphs} are finite and are permitted to have loops and parallel edges.
A graph is a {\em minor} of another if the first can be obtained from 
a subgraph of the second by contracting edges.
The cornerstone of the Graph Minors project of Robertson and Seymour is the following
excluded minor theorem. (The missing definitions are as in~\cite{rs XVI} and are
given at the end of this section.)

\begin{theorem}[{\cite[Theorem~(1.3)]{rs XVI}}]
Let $L$ be a graph.
Then there exist integers $\kappa, \rho, \xi>0$ such that every graph $G$ with no $L$-minor 
can be constructed by clique-sums, starting from graphs that are an $\le\xi$-extension of
an outgrowth by $\le\kappa$ $\rho$-rings of a graph that can be drawn in a surface in which
$L$ cannot be drawn.
\end{theorem}

In this paper we are concerned with excluding topological minors.
The first such theorem was obtained by Grohe and Marx.

\begin{theorem}[{\cite[Corollary~4.4]{gm}}] \label{gm structure}
For every graph $H$ there exist integers $\xi,\kappa,\rho,g,D$ such that
every graph $G$ with no $H$-subdivision can  be constructed by clique-sums, starting from graphs that are an $\le\xi$-extension of
either
\begin{itemize}
\item[(a)] a graph of maximum degree at most $D$, or
\item[(b)]
an outgrowth by $\le\kappa$ $\rho$-rings of a graph that can be drawn in a surface of genus at most $g$.
\end{itemize}
\end{theorem}

Thus the second outcome includes graphs drawn on surfaces in which $H$ {\em can} be drawn.
Dvo\v{r}\'ak~\cite[Theorem~3]{d} strengthened the result by restricting the graphs in (b)
to those that can be drawn in a surface $\Sigma$ in which $H$ can possibly be drawn, 
but only ``in a way in which $H$ cannot be drawn in $\Sigma$".
We omit the precise statement of Dvo\v{r}\'ak's  theorem, because it requires a large amount of
definitions that we otherwise do not need. 
Instead, let us remark that the meaning of ``the way in which $H$ cannot be drawn in $\Sigma$"
has to do with the function mf, defined as follows.

Let $H$ be a graph and $\Sigma$ a surface in which $H$ can be embedded.
We define  $\mf(H,\Sigma)$ as the minimum of $|S|$, over all embeddings of $H$ in $\Sigma$ 
and all sets $S$ of regions of the embedded graph such that
every vertex of $H$ of degree at least four is incident with a region in $S$.
When $H$ cannot be embedded in $\Sigma$, we define $\mf(H,\Sigma)$ to be infinity.

Our objective is to strengthen the theorems of Grohe and Marx, and Dvo\v{r}\'ak by reducing the value of the
constant $D$ to the maximum degree of $H$, which is best possible.
(We will prove that this value $D$ cannot be replaced by any number smaller than the maximum degree of $H$ in Section~\ref{sec:bestmaxdeg}.)
However, we are not able to extend the theorems verbatim;  our theorem gives
a structure relative to a tangle, as follows.
(Tangles, vortices and segregations are defined in Section~\ref{sec:tangles}.)

\begin{theorem} \label{main}
Let $d \geq 4$ and $h>0$ be integers.
Then there exist integers $\theta, \kappa, \rho, \xi, g \geq 0$ such that the following holds.
If $H$ is a graph of maximum degree $d$ on $h$ vertices, and a graph $G$ does not admit an $H$-subdivision, 
then for every tangle $\T$ in $G$ 
of order at least $\theta$ 
there exists a set $Z \subseteq V(G)$ with $\lvert Z \rvert \leq \xi$ such that either
		\begin{enumerate}
			\item for every vertex $v \in V(G)-Z$ there exists $(A,B) \in \T-Z$ of order at most $d-1$ such that $v \in V(A)-V(B)$, or
		\item there exists a $(\T-Z)$-central segregation $\Se = \Se_1 \cup \Se_2$ of $G-Z$ with $\lvert \Se_2 \rvert \leq \kappa$
such that $\Se$ has a proper arrangement in some surface $\Sigma$ of genus at most $g$, 
every society $(S_1,\Omega_1)$ in $\Se_1$ satisfies $\lvert \overline{\Omega}_1 \rvert \leq 3$, 
every society $(S_2, \Omega_2)$ in $\Se_2$ is a $\rho$-vortex, and either
			\begin{enumerate}
				\item $H$ cannot be drawn in $\Sigma$, or
				\item $H$ can be drawn in $\Sigma$ and $\mf(H,\Sigma) \geq 2$, and there exists $\Se'_2 \subseteq \Se_2$ with $\lvert \Se'_2 \rvert \leq \mf(H,\Sigma)-1$ such that 
for every vertex  $v\in V(G)-Z$ 
either $v\in V(S)-\bar{\Omega}$ for some $(S,\Omega) \in \Se'_2$
or there exists $(A,B) \in \T-Z$ of order at most $d-1$ such that $v \in V(A)-V(B)$.
			\end{enumerate}
		\end{enumerate}
\end{theorem}

\noindent
In fact, we will prove a stronger statement (Theorem~\ref{stronger excluding subdivision}) that provides the structure information for graphs with no $H$-subdivision with branch vertices prescribed and immediately implies Theorem \ref{main}.
In addition, Theorem~\ref{main} has the following immediate corollary.

\begin{corollary}
\label{corol}
Let $d \geq 4$ and $h>0$ be integers.
Then there exist $\theta$ and $\xi$ such that for every graph $H$ of order $h$ and of maximum degree $d$ that can be drawn in the plane such that every vertex of degree at least four is incident with the infinite region, and for every graph $G$, either $G$ admits an $H$-subdivision, or for every tangle $\T$ of order at least $\theta$ in $G$, there exists $Z \subseteq V(G)$ with $\lvert Z \rvert \leq \xi$ such that for every vertex $v \in V(G)-Z$ there exists $(A,B) \in \T-Z$ of order at most $d-1$ such that $v \in V(A)-V(B)$.
\end{corollary}

\begin{pf}
Let $d\ge4$ and $h$ be given,  let $\theta$ and $\xi$ be as in Theorem~\ref{main}, and
let $H$ be as in the statement of the corollary.
Then $\mf(H,\Sigma)=1$ for every surface $\Sigma$, and hence the second outcome of Theorem \ref{main} cannot hold.
Thus the first outcome holds, as desired.
\end{pf}

\bigskip

Corollary~\ref{corol}  will be used in a forthcoming series of papers to prove the following theorem, conjectured by Robertson.
In the application it will be important that the order of the separation in Corollary~\ref{corol} is at most $d-1$.

\begin{theorem}[\cite{lt}]
Let $k\ge1$ be an integer,  let $R$ denote the graph obtained from a path of length $k$ by replacing each
edge by a pair of parallel edges, and let $G_1,G_2,\ldots$ be an infinite sequence of graphs such that 
none of them has an $R$-subdivision. Then there exist integers $i,j$ such that $1\le i<j$ and $G_j$ has
a $G_i$-subdivision.
\end{theorem}

Let us now introduce the missing definitions.
Given a subset $X$ of the vertex-set $V(G)$ of a graph $G$, the subgraph of $G$ induced by $X$ is denoted by $G[X]$.
We say that a graph $G$ is the {\it clique-sum} of graphs $G_1,G_2$ if there exist $V_1=\{v_{1,1},...,v_{1,\lvert V_1 \rvert}\} \subseteq V(G_1), V_2 = \{v_{2,1},v_{2,2},...,v_{2,\lvert V_2 \rvert}\} \subseteq V(G_2)$ with $\lvert V_1 \rvert = \lvert V_2 \rvert$ such that $G_1[V_1]$ and $G_2[V_2]$ are complete graphs, and $G$ can be obtained from $G_1 \cup G_2$ by identifying $v_{1,i}$ and $v_{2,i}$ for each $i$ and deleting a subset of edges with both ends in $V_1\cup V_2$.
A graph $G'$ is a {\it $\leq r$-extension} of a graph $G$ if $G$ can be obtained from $G'$ by deleting at most $r$ vertices of $G$.
A graph $G$ is an {\it $r$-ring with perimeter $t_1,...,t_n$} if $t_1,...,t_n \in V(G)$ are distinct and there is a sequence $X_1,...,X_n$ of subsets of $V(G)$ such that 
	\begin{itemize}
		\item $X_1 \cup ... \cup X_n = V(G)$, and every edge of $G$ has both ends in some $X_i$, 
		\item $t_i \in X_i$ for $1 \leq i \leq n$, 
		\item $X_i \cap X_k \subseteq X_j$ for $1 \leq i \leq j \leq k \leq n$,
		\item $\lvert X_i \rvert \leq r$ for $1 \leq i \leq n$.
	\end{itemize}
Let $G_0$ be a graph drawn in a surface $\Sigma$, and let $\Delta_1,...,\Delta_d \subseteq \Sigma$ be pairwise disjoint closed disks, each meeting the drawing only in vertices of $G_0$, and each containing no vertices of $G_0$ in its interior.
For $1 \leq i \leq d$, let the vertices of $G_0$ in the boundary of $\Delta_i$ be $t_1,...,t_n$ say, in order, and choose an $r$-ring $G_i$ with perimeter $t_1,...,t_n$ meeting $G_0$ just in $t_1,...,t_n$ and disjoint from every other $G_j$; and let $G$ be the union of $G_0,G_1,...,G_d$.
We call such a graph $G$ an {\it outgrowth by $d$ $r$-rings of $G_0$}.

The paper is organized as follows.
In Section~\ref{sec:tangles} we review the notions of tangles and graph minors.
In Section~\ref{sec:spider theorem} we prove an Erd\H{o}s-P\'osa-type result for ``spiders",
trees with one vertex of degree $d$ and all other vertices of degree one or two.
In Section~\ref{sec:taming spiders} we prove a lemma that will allow us to find a large well-behaved family
of spiders, given a huge number of spiders.
In Section~\ref{sec:surfaces} we review some theorems related to graphs embedded on a surface, and prove some other lemmas.
In Section~\ref{sec:structures} we prove a structure theorem for excluding subdivisions of a fixed graph with branch vertices prescribed (Theorem~\ref{stronger excluding subdivision}), which immediately implies Theorem~\ref{main}.
In Section~\ref{sec:bestmaxdeg} we prove that the order of the separations mentioned in the first conclusion of Theorem~\ref{main} cannot be decreased, and the constants $D$ mentioned in the first conclusion of Theorem~\ref{gm structure} and the first conclusion of \cite[Theorem 3]{d} cannot be improved to be a number less than the constant $d-1$ in the first conclusion of Theorem~\ref{main}.

\section{Tangles and minors}
\label{sec:tangles}
In this section, we review some theorems about tangles and graph minors.

A {\it separation} of a graph $G$ is a pair $(A,B)$ of subgraphs with $A \cup B=G$ and $E(A \cap B) = \emptyset$, and the {\it order} of $(A,B)$ is $\lvert V(A) \cap V(B) \rvert$.
A {\it tangle} $\T$ in $G$ of {\it order} $\theta$ is a set of separations of $G$, each of order less than $\theta$ such that
\begin{enumerate}
	\item[(T1)] for every separation $(A,B)$ of $G$ of order less than $\theta$, either $(A,B) \in \T$ or $(B,A) \in \T$;
	\item[(T2)] if $(A_1, B_1), (A_2,B_2), (A_3,B_3) \in \T$, then $A_1 \cup A_2 \cup A_3 \neq G$;
	\item[(T3)] if $(A,B) \in \T$, then $V(A) \neq V(G)$.
\end{enumerate}
The notion of tangle was first defined by Roberson and Seymour in \cite{rs X}.
(T1), (T2) and (T3) are called the first, second and third tangle axiom, respectively.
In addition, we say that $G$ {\it contains} $\T$ in this case.
Furthermore, for $Z \subseteq V(G)$ with $\lvert Z \rvert<\theta$, we define $\T-Z$ to be the set of all separations $(A',B')$ of $G-Z$ of order less than $\theta-\lvert Z \rvert$ such that there exists $(A,B) \in \T$ with $Z \subseteq V(A \cap B)$, $A'=A-Z$ and $B'=B-Z$.
It is proved in \cite[Theorem 8.5]{rs X} that $\T-Z$ is a tangle in $G-Z$ of order $\theta-Z$.

Given a graph $H$, an $H$-{\it minor} of a graph $G$ is a map $\alpha$ with domain $V(H) \cup E(H)$ such that the following hold.
\begin{itemize}
	\item $\alpha(h)$ is a nonempty connected subgraph of $G$, for every $h \in V(H)$.
	\item If $h_1$ and $h_2$ are different vertices of $H$, then $\alpha(h_1)$ and $\alpha(h_2)$ are disjoint.
	\item For each edge $e$ of $H$,
		\begin{itemize}
			\item if $e$ is not a loop, $\alpha(e)$ is an edge of $G$ with one end in $\alpha(h_1)$ and one end in $\alpha(h_2)$, where $h_1,h_2$ are the ends of $e$; 
			\item if $e$ is a loop, then $\alpha(e) \in E(G)-E(\alpha(h_1))$ and every end of $\alpha(e)$ is in $\alpha(h_1)$, where $h_1$ is the end of $e$.
		\end{itemize}
	\item If $e_1, e_2$ are two different edges of $H$, then $\alpha(e_1) \neq \alpha(e_2)$.
\end{itemize}
We say that {\it $G$ contains an $H$-minor} if such a function $\alpha$ exists.
For every $h \in V(H)$, $\alpha(h)$ is called a {\it branch set} of $\alpha$.
A tangle $\T$ in $G$ {\it controls} an $H$-minor $\alpha$ if $\alpha$ is an $H$-minor such that there does not exist $(A,B) \in \T$ of order less than $\lvert V(H) \rvert$ and $h \in V(H)$ such that $V(\alpha(h)) \subseteq V(A)$.

The following theorem offers a way to obtain a tangle in a graph from a minor.

\begin{theorem}[{\cite[Theorem (6.1)]{rs X}}] \label{tangle induced by minor}
Let $G$ and $H$ be graphs.
Let $\T'$ be a tangle in $H$ of order $\theta \geq 2$.
If $G$ admits an $H$-minor $\alpha$, and $\T$ is the set of separations $(A,B)$ of $G$ of order less than $\theta$ such that there exists $(A',B') \in \T'$ with $\alpha(E(A')) = E(A) \cap \alpha(E(H))$, then $\T$ is a tangle in $G$ of order $\theta$.
\end{theorem}

The tangle $\T$ in Theorem \ref{tangle induced by minor} is called the {\it tangle induced by $\T'$}.
We say that $\T'$ is {\it conformal} with a tangle $\T''$ in $G$ if $\T \subseteq \T''$.

A {\it society} is a pair $(S,\Omega)$, where $S$ is a graph and $\Omega$ is a cyclic permutation of a subset $\bar{\Omega}$ of $V(S)$.
Let $\rho$ be a nonnegative integer.
A society $(S,\Omega)$ is a {\it $\rho$-vortex} if for all distinct $u,v \in \bar{\Omega}$, there do not exist $\rho+1$ mutually disjoint paths of $S$ between $I \cup \{u\}$ and $J \cup \{v\}$, where $I$ is the set of vertices in $\bar{\Omega}$ after $u$ and before $v$ in $\Omega$, and $J$ is the set of vertices in $\bar{\Omega}$ after $v$ and before $u$ in $\Omega$.

A {\it segregation} of a graph $G$ is a set $\Se$ of societies such that the following hold. 
\begin{itemize}
	\item $S$ is a subgraph of $G$ for every $(S, \Omega) \in \Se$, and $\bigcup \{S: (S,\Omega) \in \Se\}=G$.
	\item For every distinct $(S,\Omega)$ and $(S', \Omega') \in \Se$, $V(S \cap S') \subseteq \bar{\Omega} \cap \overline{\Omega}'$ and $E(S \cap S') = \emptyset$.
\end{itemize}
We write $V(\Se) = \bigcup \{\bar{\Omega}: (S, \Omega) \in \Se\}$.
If $\T$ is a tangle in $G$, a segregation $\Se$ of $G$ is {\it $\T$-central} if for every $(S,\Omega) \in \Se$, there is no $(A,B) \in \T$ of order at most half of the order of $\T$ with $B \subseteq S$.
Note that our definition of a $\T$-central segregation is different from the one in \cite{rs XVI}, but every segregation that is $\T$-central in the sense of \cite{rs XVI} is $\T$-central in the sense of this paper.

A {\it surface} is a nonnull compact connected $2$-manifold without boundary.
Let $\Sigma$ be a surface and $\Se = \{(S_1, \Omega_1), ..., (S_k, \Omega_k)\}$ a segregation of $G$.
For every subset $\Delta$ of $\Sigma$, we denote the closure of $\Delta$ by $\bar{\Delta}$, and the boundary of $\Delta$ by $\partial\Delta$.
An {\it arrangement} of $\Se$ in $\Sigma$ is a function $\alpha$ with domain $\Se \cup V(\Se)$, such that the following hold.
\begin{itemize}
	\item For $1 \leq i \leq k$, $\alpha(S_i, \Omega_i)$ is a closed disk $\Delta_i \subseteq \Sigma$, and $\alpha(x) \in \partial\Delta_i$ for each $x \in \overline{\Omega_i}$.
	\item For $1 \leq i \leq k$, if $x \in \Delta_i \cap \Delta_j$, then $x=\alpha(v)$ for some $v \in \overline{\Omega_i} \cap \overline{\Omega_j}$.
	\item For all distinct $x,y \in V(\Se)$, $\alpha(x) \neq \alpha(y)$.
	\item For $1 \leq i \leq k$, $\Omega_i$ is mapped by $\alpha$ to a natural order of $\alpha(\overline{\Omega_i})$ determined by $\partial\Delta_i$.
\end{itemize}
An arrangement is {\it proper} if $\Delta_i \cap \Delta_j = \emptyset$ for all $1 \leq i < j \leq k$ such that $\lvert \overline{\Omega_i} \rvert, \lvert \overline{\Omega_j} \rvert >3$.

Given a graph $H$, an {\it $H$-subdivision} is a pair of functions $(\pi_V, \pi_E)$ such that the following hold.
\begin{itemize}
	\item $\pi_V: V(H) \rightarrow V(G)$ is an injective function.
	\item $\pi_E$ maps loops of $H$ to cycles in $G$ and maps other edges of $H$ to paths in $G$ such that $\pi_E(e)$ contains $\pi_V(v)$, and $\pi_E(e')$ has ends $\pi_V(x)$ and $\pi_V(y)$ for every loop $e$ with end $v$ and every edge $e' \in E(H)$ with distinct ends $x$ and $y$.
	\item If $f_1, f_2$ are two different edges in $H$, then $\pi_E(f_1) \cap \pi_E(f_2) \subseteq \pi_V(X)$, where $X$ is the set of common ends of $f_1$ and $f_2$. 
\end{itemize}
We say that {\it $G$ admits an $H$-subdivision} if such a pair of functions $(\pi_V, \pi_E)$ exists.
The vertices in the image of $\pi_V$ are called the {\it branch vertices} of $(\pi_V,\pi_E)$.

\section{Finding disjoint spiders}
\label{sec:spider theorem}

First, we introduce a lemma proved by Robertson and Seymour \cite{rs XIII}.

\begin{lemma}[{\cite[Theorem (5.4)]{rs XIII}}] \label{vertex linkage}
Let $G$ be a graph, and let $Z$ be a subset of $V(G)$ with $\lvert Z \rvert = \xi$.
Let $k \geq \lfloor \frac{3}{2} \xi \rfloor$, and let $\alpha$ be a $K_k$-minor in $G$.
If there is no separation $(A,B)$ of $G$ of order less than $\lvert Z \rvert$ such that $Z \subseteq V(A)$ and $A \cap \alpha(h)=\emptyset$ for some $h \in V(K_k)$, then for every partition $(Z_1,...,Z_n)$ of $Z$ into non-empty subsets, there are $n$ connected graphs $T_1, ..., T_n$ of $G$, mutually disjoint and $V(T_i) \cap Z = Z_i$ for $1 \leq i \leq n$.
\end{lemma}

A {\it $d$-spider with head $v$} is a tree such that every vertex other than $v$ in the tree has degree at most $2$, and the degree of $v$ is $d$.
A {\it leaf} is a vertex of degree one.
Let $G$ be a graph, and let $S,Y$ be subsets of $V(G)$.
A {\it $d$-spider from $S$ to $Y$} is a $d$-spider with head $v \in S$ whose leaves are in $Y$.

Let $G$ be a graph and $\T$ a tangle in $G$.
We say that a subset $X$ of $V(G)$ is {\it free} if there exists no $(A,B) \in \T$ of order less than $\lvert X \rvert$ such that $X \subseteq V(A)$.

\begin{lemma}  \label{vertex spider subdivision}
Let $G$ be a graph and $H$ be a graph on $h$ vertices of maximum degree at most $d$.
Let $t \geq \lfloor \frac{3 hd}{2}\rfloor$.
Let $\T$ be a tangle of order at least $hd$ in $G$ that controls a $K_t$-minor.
If there exist pairwise disjoint sets $X_1,X_2,...,X_h$ such that for $1 \leq i \leq h$ the set $X_i$ consists of a vertex $v_i$ of $G$ and $d-1$ of its neighbors and $\bigcup_{i=1}^h X_i$ is free with respect to $\T$, then $G$ has an $H$-subdivision with branch vertices $v_1,v_2,...,v_h$.
\end{lemma}

\begin{pf}
Let $Z= \bigcup_{i=1}^h X_i$, and let $\alpha$ be a $K_t$-minor controlled by $\T$.
Suppose that there exists a separation $(A,B)$ of $G$ of order less than $\lvert Z \rvert$ such that $Z \subseteq V(A)$ and $A \cap \alpha(v) = \emptyset$ for some $v \in V(K_t)$.
By the first tangle axiom, either $(A,B) \in \T$ or $(B,A) \in \T$.
Since $Z$ is free, $(B,A) \in \T$.
But it is a contradiction since $t \geq hd$ and $\T$ controls $\alpha$.
Therefore, there does not exist a separation $(A,B)$ of $G$ of order less than $\lvert Z \rvert$ such that $Z \subseteq V(A)$ and $A \cap \alpha(v) = \emptyset$ for some $v \in V(K_t)$.

Denote $V(H)$ by $\{u_1,u_2,...,u_h\}$ and $E(H)$ by $\{e_1,e_2,...,e_{\lvert E(H) \rvert}\}$.
Since the maximum degree of $H$ is at most $d$, there exist $Z_0 \subseteq Z$ and a partition $(Z_1,Z_2,...,Z_{\lvert E(H) \rvert})$ of $Z-Z_0$ such that for every $1 \leq \ell \leq \lvert E(H) \rvert$, $Z_\ell$ consists of two distinct vertices where one is in $X_i$ and one is in $X_j$, where the ends of $e_\ell$ are $u_i$ and $u_j$.
By Lemma \ref{vertex linkage}, there exist $\lvert E(H) \rvert$ pairwise disjoint paths in $G-Z_0$ connecting the two vertices of each part of $(Z_1,Z_2,...,Z_{\lvert E(H) \rvert})$.
This creates a subdivision of $H$ with branch vertices $v_1,v_2,...,v_h$.
\end{pf}

\begin{theorem}[{\cite[Theorem 3.3]{mw}}] \label{mw spider}
Let $G$ be a graph and $\T$ a tangle in $G$ of order $\theta$.
Let $\{X_j \subseteq V(G): j \in J\}$ be a family of subsets of $V(G)$ indexed by $J$.
Let $d,k$ be integers with $\theta \geq (k+d)^{d+1}+d$.
If $\lvert X_j \rvert=d$ for every $j \in J$, then there exists a set $J' \subseteq J$ satisfying the following.
	\begin{enumerate}
		\item For distinct elements $j, j'$ of $J'$, $X_j$ and $X_{j'}$ are disjoint.
		\item $\bigcup_{j \in J'} X_j$ is free.
		\item If $\lvert \bigcup_{j \in J'} X_j \rvert < k$, then there exists a set $Z$ with $\bigcup_{j \in J'} X_j \subseteq Z$ and $\lvert Z \rvert \leq (k+d)^{d+1}$ satisfying that for all $j \in J$, either $X_j \cap Z \neq \emptyset$, or $X_j$ is not free in $\T-Z$.
	\end{enumerate}
\end{theorem}

\begin{theorem} \label{spider tangle}
Let $h$ and $d$ be positive integers.
Let $G$ be a graph, and let $S$ be a subset of vertices of degree at least $d-1$ in $G$.
Let $\T$ be a tangle in $G$ of order $\theta$.
If $\theta \geq (hd+1)^{d+1}+d$, then either 
	\begin{enumerate}
		\item there exist $h$ vertices $v_1,v_2,...,v_h \in S$ and $h$ pairwise disjoint subsets $X_1,X_2,...,X_h$ of $V(G)$, where $X_i$ consists of $v_i$ and $d-1$ neighbors of $v_i$ for each $1 \leq i \leq h$, such that $\bigcup_{i=1}^h X_i$ is free in $\T$, or
		\item there exists a set $C \subseteq V(G)$ with $\lvert C \rvert \leq (hd+1)^{d+1}$ such that for every $v \in S-C$, there exists $(A,B) \in \T-C$ of order less than $d$ such that $v \in V(A)-V(B)$.
	\end{enumerate}
\end{theorem}

\begin{pf}
Let $\{X_j: j \in J\}$ be the collection of the $d$-element subsets consisting of one vertex $v_j$ in $S$ and $d-1$ of its neighbors.
Applying Theorem \ref{mw spider} by further taking $k=(h-1)d+1$, we know there exists $J' \subseteq J$ such that $X_j \cap X_{j'} = \emptyset$ for every distinct $j,j'$ in $J'$, and $\bigcup_{j \in J'}X_j$ is free.
Furthermore, if $\lvert \bigcup_{j \in J'}X_j \rvert \leq (h-1)d$, there exists a set $C$ with $\bigcup_{j \in J'}X_j \subseteq C \subseteq V(G)$ and $\lvert C \rvert \leq (hd+1)^{d+1}$ satisfying that for all $j \in J$, either $X_j \cap C \neq \emptyset$, or $X_j$ is not free in $\T-C$.

Observe that if $\lvert \bigcup_{j \in J'}X_j \rvert > (h-1)d$, then $\lvert J' \rvert \geq h$ and the first statement holds.
So we assume that $\lvert \bigcup_{j \in J'} X_j \rvert \leq (h-1)d$, and we shall prove that the second statement of this theorem holds.
Let $v \in S-C$.
Suppose that there does not exist $(A,B) \in \T-C$ of order less than $d$ such that $v \in V(A)-V(B)$.
In particular, $v$ has at least $d$ neighbors in $G-C$.
Let $U$ be the collection of those $X_j$ that are disjoint from $C$ and consist of $v$ and $d-1$ neighbors of $v$.
For every member $X_j$ of $U$, we define the {\it rank} of $X_j$ to be the minimum order of a separation $(A,B) \in \T-C$ such that $X_j \subseteq V(A)$.
As no member of $U$ is free, the rank of each member of $U$ is at most $d-1$. 
Let $r$ be the maximum rank of any member of $U$, and let $X$ be a member of $U$ of rank $r$.
Let $(A,B) \in \T-C$ of order $r$ such that $X \subseteq V(A)$, and subject to that, $\lvert V(B)-V(A) \rvert$ is as small as possible.
By the assumption, $v \in V(A) \cap V(B)$ and $r \leq d-1$.
On the other hand, there exist $r$ disjoint paths in $G-C$ from $X-\{v\}$ to $V(B)$, as $v$ is adjacent to all vertices in $X-\{v\}$.
We denote these $r$ disjoint paths by $P_1,P_2,...,P_r$, and denote the end of $P_i$ in $X-\{v\}$ by $u_i$ for $1 \leq i \leq r$.
As $v \in V(A) \cap V(B)$ and $\lvert V(A) \cap V(B) \rvert=r$, $v \in V(P_i)$ for some $1 \leq i \leq r$.
Without loss of generality, we may assume that $v \in V(P_r)$.
In addition, $v$ is adjacent to a vertex $u$ in $V(B)-V(A)$, otherwise, the rank of $X$ is smaller than $r$.
As $(X-\{u_r\}) \cup \{u\}$ is a member of $U$, its rank is at most $r$.
Let $(A',B') \in \T-C$ be a separation of order at most $r$ such that $(X-\{u_r\}) \cup \{u\} \subseteq V(A')$.
$X \subseteq V(A \cup A')$ and $u \in (V(B)-V(A))-(V(B \cap B')-V(A \cup A'))$, so the order of $(A \cup A', B \cap B')$ is at least $r+1$ by the choice of $(A,B)$.
It implies that the order of $(A \cap A',B \cup B')$ is at most $r-1$.
Notice that $v \in V(A') \cap V(B')$ by the assumption, so $((A \cap A')-\{v\}, (B \cup B')-\{v\})$ is a separation of $G-(C \cup \{v\})$ of order less than $r-1$.
But $P_1,P_2,...,P_{r-1}$ are $r-1$ disjoint paths from $V(A \cap A')-\{v\}$ to $V(B \cup B')-\{v\}$ in $G-(C \cup \{v\})$, a contradiction.
This proves the second statement.
\end{pf}

\bigskip

We need the following variation of Theorem \ref{spider tangle}.
A version for edge-disjoint spiders was proved in \cite{m} and \cite[Theorem 3.1]{mw}.

\begin{theorem} \label{vertex spider set}
Let $G$ be a graph, and let $X$, $Y$ be disjoint subsets of $V(G)$.
Let $h,d$ be nonnegative integers.
Then either there exist $h$ disjoint $d$-spiders from $X$ to $Y$, or there exists $C \subseteq V(G)$ with $\lvert C \rvert \leq \frac{3}{2}(hd+1)^{d+1} + \frac{d}{2}+1$ such that every $d$-spider from $X$ to $Y$ intersects $C$.
\end{theorem}

\begin{pf}
Note that for every subset $C$ of $Y$ such that $\lvert Y-C \rvert \leq d-1$, every $d$-spider from $X$ to $Y$ intersects $C$.
So we may assume that $\lvert Y \rvert \geq \frac{3}{2}((hd+1)^{d+1}+d)$, otherwise we are done.
Let $G'$ be the graph obtained from $G$ by adding edges such that $Y$ is a clique in $G'[Y]$.
As every complete graph on $k$ vertices contains a tangle of order $\lfloor 2k/3 \rfloor$, $G'[Y]$ contains a tangle of order $(hd+1)^{d+1}+d$.
And $G'[Y]$ is a minor of $G'$, so $G'$ contains a tangle $\T$ of order $(hd+1)^{d+1}+d$ induced by a tangle of the same order in $G'[Y]$ by Theorem \ref{tangle induced by minor} such that $Y \subseteq V(B)$ for every $(A,B) \in \T$.
Let $\{X_j: j \in J\}$ be the collection of $d$-element subsets of $V(G)$ such that every $X_j$ consists of one vertex $x$ in $X$ and $d-1$ neighbors of $x$.
By Theorem \ref{mw spider}, there exists $J' \subseteq J$ such that $X_j \cap X_{j'} = \emptyset$ for every distinct $j,j'$ in $J'$, and $\bigcup_{j \in J'}X_j$ is free.
Furthermore, if $\lvert \bigcup_{j \in J'}X_j \rvert \leq (h-1)d$, there exists $C \subseteq V(G)$ with $\lvert C \rvert \leq (hd+1)^{d+1}$ satisfying that for all $j \in J$, either $X_j \cap C \neq \emptyset$, or $X_j$ is not free in $\T-C$.

First, assume that $\lvert \bigcup_{j \in J'}X_j \rvert > (h-1)d$, so $\lvert J' \rvert \geq h$.
Let $\{1,2,...,h\} \subseteq J'$, and  for $1 \leq j \leq h$, let $x_j$ be a vertex in $X_j \cap X$ adjacent to all other vertices in $X_j$.
Suppose that there do not exist $dh$ disjoint paths from $\bigcup_{j=1}^h X_j$ to $Y$ in $G'$.
Then there exists a separation $(A,B)$ of $G'$ of order less than $dh$ such that $\bigcup_{j=1}^h X_j \subseteq V(A)$ and $Y \subseteq V(B)$.
Since $Y \subseteq V(B)$, we know that $(A,B) \in \T$.
But it implies that $\bigcup_{j=1}^h X_j$ is not free, a contradiction.
Hence, there exist $dh$ disjoint paths from $\bigcup_{j=1}^h X_j$ to $Y$ in $G'$.
That is, there exist $h$ disjoint $d$-spiders from $\{x_j: 1 \leq j \leq h\}$ to $Y$ in $G'$.
We are done in this case since every $d$-spider from $X$ to $Y$ in $G'$ contains a $d$-spider from $X$ to $Y$ in $G$ as a subgraph.

So we may assume that $\lvert \bigcup_{j \in J'}X_j \rvert \leq (h-1)d$, there exists $C \subseteq V(G)$ with $\bigcup_{j \in J'}X_j \subseteq C$ and $\lvert C \rvert \leq (hd+1)^{d+1}$ satisfying that for all $j \in J$, either $X_j \cap C \neq \emptyset$, or $X_j$ is not free in $\T-C$.
Let $v \in X-C$, and let $D$ be a $d$-spider from $v$ to $Y$ in $G$.
Note that $D$ is also a $d$-spider from $v$ to $Y$ in $G'$.
Suppose that $D$ is disjoint from $C$.
So $V(D)$ contains some $X_j$ such that $v \in X_j$ and $X_j \cap C = \emptyset$.
Since $X_j$ is not free in $\T-C$, there exists $(A,B) \in \T-C$ of order less than $d$ such that $X_j \subseteq V(A)$ and $Y-C \subseteq V(B)$.
It is a contradiction since there exist $d$ disjoint paths in $D$ from $V(A)$ to $V(B)$.
This proves that $D$ intersects $C$.
\end{pf}

\section{Taming spiders}
\label{sec:taming spiders}

A {\it surface} is a compact $2$-manifold.
An {\it O-arc} is a subset homeomorphic to a circle, and a {\it line} is a subset homeomorphic to $[0,1]$.
Let $\Sigma$ be a surface.
A {\it drawing} $\Gamma$ in $\Sigma$ is a pair $(U,V)$, where $V \subseteq U \subseteq \Sigma$, $U$ is closed, $V$ is finite, $U-V$ has only finitely many arc-wise connected components, called {\it edges}, and for every edge $e$, either $\bar{e}$ is a line whose set of ends are $\bar{e} \cap V$, or $\bar{e}$ is an O-arc and $\lvert \bar{e} \cap V \rvert =1$.
The components of $\Sigma-U$ are called {\it regions}.
The members of $V$ are called {\it vertices}.
For a drawing $\Gamma=(U,V)$, we write $U=U(\Gamma), V= V(\Gamma)$, and $E(\Gamma),R(\Gamma)$ are defined to be the set of edges and the set of regions, respectively. 
If $v$ is a vertex of a drawing $\Gamma$ and $e$ is an edge or a region of $\Gamma$, we say that $e$ is {\it incident with} $v$ if $v$ is contained in the closure of $e$.
Note that the incidence relation between $V(\Gamma)$ and $E(\Gamma)$ defines a graph, and we say that $\Gamma$ is a {\it drawing of $G$} in $\Sigma$ if $G$ is defined by this incidence relation.
In this case, we say that $G$ is {\it embeddable} in $\Sigma$, or $G$ can be {\it drawn} in $\Sigma$.

We say that $(S,\Omega, \Omega_0)$ is a {\it neighborhood} if $S$ is a graph and $\Omega, \Omega_0$ are cyclic permutations with $\bar{\Omega}, \overline{\Omega_0} \subseteq V(S)$.
A neighborhood $(S,\Omega, \Omega_0)$ is {\it rural} if $S$ has a drawing $\Gamma$ in the plane and there are disks $\Delta_0 \subseteq \Delta$ such that 
	\begin{itemize}
		\item $\Gamma$ uses no point outside $\Delta$ and none in the interior of $\Delta_0$, and 
		\item $\bar{\Omega}$ are the vertices in $\Gamma \cap \partial\Delta$, and $\overline{\Omega_0}$ are the vertices in $\Gamma \cap \Delta_0$, and 
		\item the cyclic permutations of $\bar{\Omega}$ and $\overline{\Omega_0}$ coincide with the natural cyclic order on $\Delta$ and $\Delta_0$.
	\end{itemize}
In this case, we say that $(\Gamma, \Delta, \Delta_0)$ is a {\it presentation} of $(S,\Omega, \Omega_0)$.
For a fixed presentation $(\Gamma, \Delta, \Delta_0)$ of a neighborhood $(S,\Omega, \Omega_0)$ and an integer $s \geq 0$, an {\it $s$-nest} for $(\Gamma, \Delta, \Delta_0)$ is a sequence $(C_1, C_2, ..., C_s)$ of pairwise disjoint cycles of $S$ such that $\Delta_0 \subseteq \Delta_1 \subseteq ... \subseteq \Delta_s \subseteq \Delta$, where $\Delta_i$ is the closed disk in the plane bounded by $C_i$ in the drawing $\Gamma$.

If $(S,\Omega, \Omega_0)$ is a neighborhood and $(S_0,\Omega_0)$ is a society, then $(S \cup S_0, \Omega)$ is a society and we call this society the {\it composition} of the society $(S_0,\Omega_0)$ with the neighborhood $(S,\Omega, \Omega_0)$.
A society $(S,\Omega)$ is {\it $s$-nested} if it is the composition of a society with a rural neighborhood that has an $s$-nest for some presentation of it.

A subgraph $F \subseteq S$ for a rural neighborhood $(S,\Omega,\Omega_0)$ with presentation $(\Gamma,\Delta,\Delta_0)$ is {\it perpendicular} to an $s$-nest $(C_1,C_2,...,C_s)$ for $(\Gamma,\Delta,\Delta_0)$ if for every component $P$ of $F$
	\begin{itemize}
		\item $P$ is a path with one end in $\bar{\Omega}$ and the other in $\overline{\Omega_0}$, and 
		\item $P \cap C_i$ is a path for all $i=1,2,...,s$.
	\end{itemize}

We shall use the following theorem, which was proved in \cite{kntw}, to prove the main theorem of this section.
We present a simplified restatement of it.

\begin{theorem}[{\cite[Theorem 10.5]{kntw}}] \label{buffer and target}
For every three positive integers $s,k \geq 3,c$, there exists an integer $s'(s,k,c)$ such that for every $s'$-nested society $(S,\Omega)$ that is a composition of a society $(S_0,\Omega_0)$ with a rural neighborhood with an $s'$-nest, and for every union of $c$ pairwise disjoint $k$-spiders $F_0$ from $V(S_0)-\overline{\Omega_0}$ to $\bar{\Omega}$, where every vertex in $V(F_0) \cap \bar{\Omega}$ is a leaf in $F_0$, there exists a union of $c$ pairwise disjoint $k$-spiders $F$ in $(S,\Omega)$ from the set of the heads of $F_0$ to the set of leaves of $F_0$ such that $(S,\Omega)$ can be expressed as a composition of some society with a rural neighborhood $(S',\Omega,\Omega')$ that has a presentation with an $s$-nest $(C_1,C_2,...,C_s)$ such that $S' \cap F$ is perpendicular to $(C_1,C_2,...,C_s)$.
\end{theorem}

Given a cyclic ordering $\Omega$ on a set $\bar{\Omega}$, an {\it interval} of $\Omega$ is a subset $I$ of $\bar{\Omega}$ such that the elements of $I$ are consecutive elements in $\Omega$.

A {\it tree-decomposition} of a graph $G$ is a pair $(T,\X)$, where $T$ is a tree and $\X$ is a collection $\{X_t: t \in V(G)\}$ of subsets of $V(G)$, such that the following hold.
	\begin{itemize}
		\item $\bigcup_{t \in V(T)} X_t = V(G)$.
		\item For each edge $e$ of $G$, there exists $t \in V(T)$ such that $X_t$ contains the ends of $e$.
		\item For each vertex $v$ of $G$, the subgraph of $T$ induced by $\{t: v \in X_t\}$ is connected.
	\end{itemize}
The {\it width} of $(T,\X)$ is $\max\{\lvert X_t \rvert-1: t \in V(T)\}$.
The {\it adhesion} of $(T,\X)$ is $\max\{\lvert X_t \cap X_{t'} \rvert: tt' \in E(T)\}$.
For each $t \in V(T)$, the {\it torso} at $t$ is the graph $L$ obtained from the subgraph of $G$ induced by $X_t$ by adding edges such that for each neighbor $t'$ of $t$ in $T$, $X_t \cap X_{t'}$ is a clique in $L$.
When the tree $T$ is a path, we say $(T,\X)$ is a {\it path-decomposition}.
The {\it tree-width} of $G$ is the minimum width of a tree-decomposition of $G$.

Now, we are ready to prove the main theorem of this section.

\begin{theorem} \label{nest buffer makes good}
For every positive integers $d \geq 3, \rho,k$ and $s$, there exist integers $s'=s'(k,d,s,\rho)$ and $k'=k'(k,\rho)$ such that for every $s'$-nested society $(S,\Omega)$ that is a composition of a $\rho$-vortex $(S_0,\Omega_0)$ with a rural neighborhood that has an $s'$-nest, and for every $k'$ pairwise disjoint $d$-spiders $D_1,D_2,...,D_{k'}$ from $V(S_0)-\overline{\Omega_0}$ to $\bar{\Omega}$ such that every vertex of $D_i \cap \bar{\Omega}$ is a leaf of $D_i$, there exist $k$ pairwise disjoint $d$-spiders $D'_1, D'_2,...,D'_k$ from $V(S_0)$ to $\bar{\Omega}$ such that the following hold.
	\begin{enumerate}
		\item  $(S,\Omega)$ can be expressed as a composition of a society $(S_0',\Omega')$ with a rural neighborhood $(S',\Omega,\Omega')$ that has a presentation with an $s$-nest $(C_1,C_2,...,C_s)$ such that $D_i' \cap S'$ is perpendicular to $(C_1,C_2,...,C_s)$ for every $1 \leq i \leq k$.
		\item For every $1 \leq i \leq k$, the head of $D'_i$ is the head of $D_{i'}$ for some $1 \leq i' \leq k'$.
		\item For every $1 \leq i \leq k$, every leaf of $D_i'$ is a leaf of $D_1 \cup D_2 \cup ... \cup D_{k'}$.
		\item For every $1 \leq i \leq k$, there exists an interval $I_i$ of $\Omega$ containing all leaves of $D_i'$ such that $I_i$ is disjoint from $I_j$ for every $j \neq i$.
	\end{enumerate}
\end{theorem}

\begin{pf}
Let $s'(k,d,s,\rho) = s'_{\ref{buffer and target}}(s,d,(12\rho+7)k)$ and $k'(k,\rho) = (12\rho+7)k$, where $s'_{\ref{buffer and target}}$ is the function $s'$ mentioned in Theorem \ref{buffer and target}.
By Theorem \ref{buffer and target}, there exist $(12\rho+7)k$ pairwise disjoint $d$-spiders $D'_1,D'_2,...,D'_{k'}$ from the set of the heads of $D_1,D_2,...,D_{k'}$ to the union of the set of leaves of $D_1,D_2,...,D_{k'}$ such that $(S,\Omega)$ can be expressed as a composition of some society with a rural neighborhood $(S',\Omega,\Omega')$ that has a presentation with an $s$-nest $(C_1,C_2,...,C_s)$ such that $D_i' \cap S'$ is perpendicular to $(C_1,C_2,...,C_s)$ for every $1 \leq i \leq k'$.

Let $\overline{\Omega_0} = \{v_1,v_2,...,v_{\lvert \overline{\Omega_0} \rvert}\}$ in order.
Since the head of each $D_i'$ is contained in $V(S_0)-\overline{\Omega_0}$, each component of $D_i'-(V(S_0)-\overline{\Omega_0})$ is a path.
For each $1 \leq i \leq k'$, let $W_i$ be the subset of $[\lvert \overline{\Omega_0} \rvert]$ such that for each $j \in W_i$, some component of $D_i'-(V(S_0)-\overline{\Omega_0})$ is a path from a leaf of $D_i'$ to $v_j$, and let $a_i=\min W_i$ and $b_i = \max W_i$.
For each $1 \leq i \leq k'$, define $\ell_i$ and $r_i$ to be the leaves of $D_i'$ such that there exists a path in $S-(V(S_0)-\overline{\Omega_0}) \cap D_i'$ from $\ell_i$ to $v_{a_i}$ and from $r_i$ to $v_{b_i}$, respectively; define $I_i$ to be the interval of $\Omega$ with ends $\ell_i$ and $r_i$ containing all leaves of $D_i'$.
Then it is sufficient to prove that there exist $1 \leq i_1 < i_2 <... <i_k \leq k'$ such that $I_{i_1},I_{i_2},...,I_{i_k}$ are pairwise disjoint.

Since $(S_0,\Omega_0)$ is a $\rho$-vortex, by Theorem 8.1 in \cite{rs IX}, there exists a path-decomposition $(t_1t_2...t_{\lvert \overline{\Omega_0} \rvert}, \X)$ of $S_0$, where $\X=(X_{t_i}: 1 \leq i \leq \lvert \overline{\Omega_0} \rvert)$, such that $\lvert X_{t_i} \cap X_{t_j} \rvert \leq \rho$ for every $1 \leq i < j \leq \lvert \overline{\Omega_0} \rvert$ and $v_i \in X_{t_i}$ for every $1 \leq i \leq \lvert \overline{\Omega_0} \rvert$.
And we define $X_0=\emptyset$ and $X_{\lvert \overline{\Omega_0} \rvert+1}=\emptyset$.
Furthermore, for each $1 \leq i \leq k'$, let $Q_i$ be the path in $D_i'$ from $r_i$ to $\ell_i$, $L_i$ be the path in $D_i'$ from $\ell_i$ to $v_{a_i}$, and $R_i$ be the path in $D_i'$ from $r_i$ to $v_{b_i}$.
We say that $I_i$ is {\it naughty} if $Q_i-((S_0-\overline{\Omega_0}) \cup V(L_i) \cup V(R_i))$ contains a path from $v_c$ to $v_{c'}$, where $1 \leq c <a_i<b_i <c' \leq \lvert \overline{\Omega_0} \rvert$ such that $Q_i$ passes through $\ell_i,v_c,v_{c'},r_i$ in the order listed; otherwise we say $I_i$ is {\it nice}.

\noindent{\bf Claim 1:} If $I_i \cap I_j \neq \emptyset$ for some $i \neq j$, and either both $I_i,I_j$ are naughty or both $I_i,I_j$ are nice, then there exists $a \in [\lvert \overline{\Omega} \rvert]$ such that $X_a \cap X_{a+1} \cap V(D_i') \neq \emptyset \neq X_a \cap X_{a+1} \cap V(D_j')$.

\noindent{\bf Proof of Claim 1:}
Since $V(D_i') \cap V(D_j') = \emptyset$, we may assume that $b_i<b_j$ by symmetry.
Since $I_i \cap I_j \neq \emptyset$, $a_j<b_i<b_j$.
Since $(S,\Omega)$ is a composition of $(S_0,\Omega_0)$ with a rural neighborhood, and $V(D_i') \cap V(D_j')=\emptyset$, the vertices $\ell_j,r_i,r_j$ appear in $\Omega$ in the order listed.
Note that $Q_i$ passes through $v_{b_i}$, a vertex in $V(S_0)-\overline{\Omega_0}$, and $v_{a_i}$ in the order listed.
So $V(Q_i) \cap X_{t_{b_i}} \cap (X_{t_{b_i-1}} \cup X_{t_{b_i+1}}) \neq \emptyset$.
Similarly, $V(Q_i) \cap X_{t_{a_i}} \cap (X_{t_{a_i-1}} \cup X_{t_{a_i+1}})$, $V(Q_j) \cap X_{t_{a_j}} \cap (X_{t_{a_j-1}} \cup X_{t_{a_j+1}})$ and $V(Q_j) \cap X_{t_{b_j}} \cap (X_{t_{b_j-1}} \cup X_{t_{b_j+1}})$ are non-empty sets.
If $V(Q_i) \cap X_{t_{b_j-1}} \cap X_{t_{b_j}} \neq \emptyset \neq V(Q_i) \cap X_{t_{b_j}} \cap X_{t_{b_j+1}}$, then we are done by choosing $a=b_j-1$ or $a=b_j$.
So we may assume that $V(Q_i) \cap X_{t_{b_j-1}} \cap X_{t_{b_j}} =\emptyset$ or $V(Q_i) \cap X_{t_{b_j}} \cap X_{t_{b_j+1}} =\emptyset$.

Let $Y_r$ be one of $X_{t_{b_j-1}} \cap X_{t_{b_j}}$ and $X_{t_{b_j}} \cap X_{t_{b_j+1}}$ such that $V(Q_i) \cap Y_r$ is empty.
Similarly, there exists $Y_\ell \in \{X_{t_{a_j-1}} \cap X_{t_{a_j}}, X_{t_{a_j}} \cap X_{t_{a_j+1}}\}$ such that $V(Q_i) \cap Y_\ell=\emptyset$; otherwise we are done by taking $a=a_j-1$ or $a=a_j$.
Hence $V(L_j) \cup V(R_j) \cup Y_\ell \cup Y_r$ is disjoint from $V(D_i')$.
So $a_j<a_i<b_i<b_j$ and $I_i$ is nice.

If $V(D_i') \cap X_{t_{b_i-1}} \cap X_{t_{b_i}}=\emptyset$, then $v_{a_i}$ and the neighbor of $v_{b_i}$ in $Q_i-V(R_i)$ belong to different components of $S-(V(R_i) \cup (X_{t_{b_i-1}} \cap X_{t_{b_i}}) \cup V(R_j) \cup Y_r)$, but $Q_i-(V(R_i) \cup (X_{t_{b_i-1}} \cap X_{t_{b_i}}) \cup V(R_j) \cup Y_r)$ contains a path passing through these two vertices, a contradiction.
So $V(D_i') \cap X_{t_{b_i-1}} \cap X_{t_{b_i}} \neq \emptyset$.
Since $I_i$ is nice, $I_j$ is nice.
If $V(D_j') \cap X_{t_{b_i-1}} \cap X_{t_{b_i}}=\emptyset$, then since $I_j$ is nice, $\ell_j$ and $r_j$ belong to different components of $D_j'-(X_{t_{b_i-1}} \cup X_{t_{b_i}})$, but $Q_j-(X_{t_{b_i-1}} \cup X_{t_{b_i}})$ contains a path passing through these two vertices, a contradiction.
Therefore, $V(D_j') \cap X_{t_{b_i-1}} \cap X_{t_{b_i}} \neq \emptyset$.
This proves the claim by taking $a=b_i-1$.
$\Box$

Suppose that there do not exist such $k$ pairwise disjoint intervals among $I_1,I_2,...,I_{k'}$.
Let $H$ be the intersection graph of $I_1,I_2,...,I_{k'}$.
Then $H$ does not contain an independent set of size $k$.
We claim that $H$ contains a clique with size at least $12\rho+7$. 
Let $u_1,u_2$ be consecutive vertices in $\Omega$.
If there are at least $12\rho+7$ members of $\{I_1,I_2,...,I_{k'}\}$ containing both $u_1,u_2$, then these $12\rho+7$ members form a clique in $H$.
So we may assume that there exist at most $12\rho+6$ members of $\{I_1,I_2,...,I_{k'}\}$ containing both $u_1,u_2$.
Let $H'$ be the intersection graph of the members of $\{I_1,I_2,...,I_{k'}\}$ not containing both $u_1,u_2$.
Then $H'$ is an interval graph on at least $k'-12\rho-6$ vertices, and $H'$ is an induced subgraph of $H$.
Since $H$ has no independent set of size $k$, $H'$ does not have an independent set of size $k$.
Since interval graphs are prefect, $H'$ (and hence $H$) has a clique with size at least $(k'-12\rho-6)/(k-1) \geq 12\rho+7$.

Therefore, there exist $12\rho+7$ pairwise intersecting members of $\{I_1,I_2,...,I_{k'}\}$. 
So there exist $6\rho+4$ pairwise intersecting members of $\{I_1,I_2,...,I_{k'}\}$ such that either all of them are naughty or all of them are nice.
Let $G$ be the graph obtained from $S$ by adding edges such that vertices in $X_i \cap X_{i+1}$ are pairwise adjacent in $G$, for every $1 \leq i \leq \lvert \overline{\Omega_0} \rvert-1$.
Since there exist $6\rho+4$ pairwise intersecting members of $\{I_1,I_2,...,I_{k'}\}$ where either all of them are naughty or all of them are nice, Claim 1 implies that $G$ contains a $K_{6\rho+4}$-minor, where each branch set is $D'_i$ for some $i$.
Without loss of generality, we may assume that the branch sets of the $K_{6\rho+4}$-minor are $D'_1,D'_2,...,D'_{6\rho+4}$.

Let $G'$ be the graph obtained from $G$ by the following operations.
	\begin{itemize}
		\item Deleting vertices not in $D'_1 \cup D'_2 \cup ... \cup D'_{6\rho+4}$.
		\item Contracting each path in $D'_i -E(S_0)$ from $\overline{\Omega_0}$ to $\overline{\Omega}$ internally disjoint from $V(S_0)$ into a vertex for every $1 \leq i \leq 6\rho+4$.
		\item For every $1 \leq i \leq 6\rho+4$ and for each path in $D_i'-E(S_0)$ from $\overline{\Omega_0}$ to $\overline{\Omega_0}$ internally disjoint from $V(S_0)$, contracting all edges incident with at least one internal vertex of this path into an edge.
		\item Repeatedly contracting edges of $D'_i[X_{t_j}]$ having an end not in $\{v_j\} \cup ((X_{t_{j-1}} \cup X_{t_{j+1}}) \cap X_{t_j})$ until each remaining edge has ends in $\{v_j\} \cup ((X_{t_{j-1}} \cup X_{t_{j+1}}) \cap X_{t_j})$ for every $1 \leq i \leq 6\rho+4$ and $1 \leq j \leq \lvert \overline{\Omega_0} \rvert$.
	\end{itemize}
Note that $G'$ contains a $K_{6\rho+4}$-minor, so the tree-width of $G'$ is at least $6\rho+3$.
Observe that $V(G') \subseteq V(S_0)$.
Furthermore, $G'$ can be written as $G_0 \cup G_1$ such that $V(G_0 \cap G_1)=V(G_0) \subseteq \overline{\Omega_0}$, and $G_0$ is an outerplanar graph that can be drawn in the plane such that the vertices of $V(G_0 \cap G_1)$ are in the boundary of a region in order, and $G_1$ has a path-decomposition of width at most $2\rho$ such that each bag is a subset of $\{v_j\} \cup (X_{t_j} \cap (X_{t_{j-1}} \cup X_{t_{j+1}}))$ (for some $j$) and contains a vertex in $V(G_0 \cap G_1)$ in order.
By Lemma 8.1 in \cite{dt}, $G'$ has tree-width less than $6\rho+3$, a contradiction.
This proves the theorem.
\end{pf}

\section{Theorems on surfaces}
\label{sec:surfaces}

In this section, we recall some results about graphs embedded in surfaces.

Let $\Sigma$ be a surface and let $\Gamma$ be a drawing in $\Sigma$.
The sets $\{v\}$, for $v \in V(\Gamma)$, the edges and regions of $\Gamma$ are called the {\it atoms} of $\Gamma$.
A {\it subdrawing} $\Gamma'$ of $\Gamma$ is a drawing satisfying $V(\Gamma') \subseteq V(\Gamma)$ and $E(\Gamma') \subseteq E(\Gamma)$; we write $\Gamma' \subseteq \Gamma$ if $\Gamma'$ is a subdrawing of $\Gamma$.
A drawing is {\it $2$-cell} if every region is an open disk.

Let $\Gamma$ be a $2$-cell drawing in a surface $\Sigma$.
We say that a drawing $K$ in $\Sigma$ is a {\it radial drawing} of $\Gamma$ if it satisfies the following conditions.
	\begin{itemize}
		\item $U(\Gamma) \cap U(K) = V(\Gamma) \subseteq V(K)$.
		\item Each region $r$ of $\Gamma$ contains a unique vertex of $K$.
		\item $K$ is a drawing of a bipartite graph, and $(V(\Gamma), V(K)-V(\Gamma))$ is a bipartition of it.
		\item For every $v \in V(\Gamma)$, the edges of $K \cup \Gamma$ incident with $v$ belong alternately to $\Gamma$ and to $K$ (in their cyclic order around $v$).
	\end{itemize}

Let $\Sigma$ be a surface, and let $\Gamma$ be a drawing in $\Sigma$.
A subset $Z$ of $\Sigma$ is {\it $\Gamma$-normal} if $Z \cap U(\Gamma) \subseteq V(\Gamma)$.
If $\Sigma$ is not a sphere, we say that $\Gamma$ is {\it $\theta$-representative} if $\lvert F \cap V(\Gamma) \rvert \geq \theta$ for every non-null-homotopic $\Gamma$-normal O-arc $F$ in $\Sigma$.
If $\Delta \subseteq \Sigma$ is a closed set such that either $\bar{e} \subseteq \Delta$ or $e \cap \Delta=\emptyset$ for each $e \in E(\Gamma)$, then we define $\Gamma \cap \Delta$ to be the drawing $(U(\Gamma) \cap \Delta, V(\Gamma) \cap \Delta)$.

Let $\Sigma$ be a surface, and let $\Gamma$ be a drawing of a graph $G$ in $\Sigma$.
A {\it tangle} in $\Gamma$ and a {\it separation} of $\Gamma$ are a tangle in $G$ and a separation of $G$, respectively.
A tangle $\T$ in $\Gamma$ of order $\theta$ is said to be {\it respectful (towards $\Sigma$)} if $\Sigma$ is connected and for every $\Gamma$-normal O-arc  $F$ in $\Sigma$ with $\lvert F \cap V(\Gamma) \rvert < \theta$, there is a closed disk $\Delta \subseteq \Sigma$ with $\partial\Delta=F$ such that $(\Gamma \cap \Delta, \Gamma \cap \overline{\Sigma-\Delta}) \in \T$.
It is clear that $\Delta$ has to be unique, and we write $\Delta = \ins(F)$; the function {\it ins} is called the {\it inside function} of $\T$.
Assume that $\Gamma$ is $2$-cell, and let $K$ be a radial drawing of $\Gamma$.
If $W$ is a closed walk of $K$, we define $K|W$ to be the subdrawing of $K$ formed by the vertices and the edges in $W$.
If the length of $W$ is less than $2\theta$, then we define $\ins(W)$ to be the union of $U(K|W)$ and $\ins(C)$, taken over all cycles $C$ of $K|W$.
For every two atoms $a,b$ of $K$, define a function $m_\T(a,b)$ as follows:
	\begin{itemize}
		\item if $a=b$, then $m_\T(a,b)=0$;
		\item if $a \neq b$ and $a,b \subseteq \ins(W)$ for some closed walk $W$ of $K$ of length less than $2 \theta$, then $m_\T(a,b) = \min\frac{1}{2} \lvert E(W) \rvert$, taking over all such closed walks $W$;
		\item otherwise, $m_\T(a,b) = \theta$.
	\end{itemize}
Note that $K$ is bipartite, so $m_\T$ is integral.
In addition, for every atom $c$ of $\Gamma$, we define $a(c)$ to be an atom of $K$ such that
	\begin{itemize}
		\item $a(c)=c$ if $c \subseteq V(\Gamma)$; 
		\item $a(c)$ is the region of $K$ including $c$ if $c$ is an edge of $\Gamma$; 
		\item $a(c)=\{v\}$, where $v$ is the vertex of $K$ in $c$, if $c$ is a region of $\Gamma$.
	\end{itemize}
For every atoms $b,c$ of $\Gamma$, we define $m_\T(b,c) = m_\T(a(b),a(c))$.
If $X,Y$ are sets of atoms of $G$, then we define $m_\T(X,Y)=\min\{m_\T(x,y): x \in X, y \in Y\}$.
When one of $X$ and $Y$, say $Y$, has size one, then we denote $m_\T(X,Y)$ by $m_\T(X,y)$, where $y$ is the unique element of $Y$.

The following is a consequence of Theorem 9.1 of \cite{rs XI}.

\begin{theorem}
Let $\Sigma$ be a surface, and let $\Gamma$ be a $2$-cell drawing of a graph in $\Sigma$.
If $\T$ is a respectful tangle in $\Gamma$, then $m_\T$ is a metric on the atoms of $\Gamma$.
\end{theorem}

The following theorem is useful.

\begin{theorem}[{\cite[Theorem~(1.1)]{rs XII}}] \label{A O-arc} 
Let $\Sigma$ be a surface, and let $\Gamma$ be a $2$-cell drawing of a graph in $\Sigma$ with $E(\Gamma) \neq \emptyset$.
Let $\T$ be a respectful tangle of order $\theta$ in $\Gamma$, and let $K$ be a radial drawing of $\Gamma$.
Let $(A,B)$ be a separation of $\Gamma$ of order less than $\theta$.
Then $(A,B) \in \T$ if and only if for every edge $e$ of $A$, there exists a cycle $C$ of $K$ with $V(C) \cap V(\Gamma) \subseteq V(A) \cap V(B)$ and with $e \subseteq \ins(C)$.
\end{theorem}

\begin{theorem} \label{A distance} 
Let $\Sigma$ be a surface, and let $\Gamma$ be a $2$-cell drawing of a graph in $\Sigma$ with $E(\Gamma) \neq \emptyset$.
Let $\T$ be a respectful tangle of order $\theta$ in $\Gamma$.
Then the following hold.
	\begin{enumerate}
		\item If $x \in V(\Gamma)$ and $(A,B) \in \T$ is a separation of $\Gamma$ such that $x \in V(A)-V(B)$ and there exists a path $P$ in $A$ from $x$ to a vertex $y \in V(A)$ internally disjoint from $V(B)$, then $m_\T(x,y) \leq \lvert V(A) \cap V(B) \rvert$.
		\item If $X \subseteq V(\Gamma)$ and $(A,B) \in \T$ is a separation of $\Gamma$ of order less than $\lvert X \rvert$ such that $X \subseteq V(A)$, subject to this, $A$ is minimal, then $m_\T(X,y) \leq \lvert V(A) \cap V(B) \rvert$ for every $y \in V(A)$.
	\end{enumerate}
\end{theorem}

\begin{pf}
We first prove Statement 1.
Let $x$ and $(A,B)$ be a vertex and a separation mentioned in the statement of Conclusion 1.
Let $y \in V(A)$ be a vertex different from $x$.
Let $P$ be a path in $A$ from $x$ to $y$ internally disjoint from $V(B)$.
Let $e$ be the edge in $P$ incident with $x$.
By Theorem \ref{A O-arc}, there exists a cycle $C$ of the radial drawing $K$ of $\Gamma$ with $V(C) \cap V(\Gamma) \subseteq V(A) \cap V(B)$ and with $e \subseteq \ins(C)$.
So $x \in \ins(C)$.
If $y \not \in \ins(C)$, then $C$ intersects $P$ at an internal vertex of $P$.
However, $V(C) \cap V(\Gamma) \subseteq V(A) \cap V(B)$.
This implies that some internal vertex of $P$ is in $V(A) \cap V(B)$, a contradiction.
Hence, $y \in \ins(C)$.
Therefore, $m_\T(x,y) \leq \lvert V(A) \cap V(B) \rvert$.

Now we prove Statement 2.
Let $X$ and $(A,B)$ be a set and a separation mentioned in the statement of Conclusion 2.
Since $\lvert V(A \cap B) \rvert < \lvert X \rvert$ and $X \subseteq V(A)$, $X-V(B) \neq \emptyset$.
Since $A$ is minimal, every component of $A-V(A \cap B)$ intersects $X-V(B)$, and each vertex in $V(A \cap B)$ is either in $X$ or adjacent to some vertex in $V(A)-V(B)$.
So for every $y \in V(A)$, either $y \in X$, or there exists a path in $A$ from a vertex $x_y \in X-V(B)$ to $y$ internally disjoint from $V(B)$.
If $y \in X$, then $m_\T(X,y)=0$; otherwise, $m_\T(x_y,y) \leq \lvert V(A \cap B) \rvert$ by Statement 1.
Hence $m_\T(X,y) \leq \lvert V(A \cap B) \rvert$ for every $y \in V(A)$.
\end{pf}

\begin{theorem}[{\cite[Theorem (8.12)]{rs XI}}, {\cite[Theorem (1.2)]{rs XII}}] \label{something is far}
Let $\T$ be a respectful tangle of order $\theta$, where $\theta \geq 2$, in a $2$-cell drawing $\Gamma$ in a connected surface $\Sigma$.
If $c$ is an atom in $\Gamma$, then there exists an edge $e$ of $\Gamma$ such that $m_\T(c,e) = \theta$.
\end{theorem}

Let $\Gamma$ be a $2$-cell drawing in a surface $\Sigma$, and let $\T$ be a respectful tangle of order $\theta$ in $\Gamma$.
Let $x$ be an atom of $\Gamma$.
A {\it $\lambda$-zone around $x$} is an open disk $\Delta$ in $\Sigma$ with $x \subseteq \Delta$, such that $\partial\Delta$ is an O-arc, $\partial\Delta \subseteq U(\Gamma)$, $m_\T(x,y) \leq \lambda$ for every atom $y$ of $\Gamma$ with $y \subseteq \bar{\Delta}$, and if $x \in E(\Gamma)$, then $\lambda \geq 2$.
A {\it $\lambda$-zone} is a $\lambda$-zone around some atom.

Let $\Delta$ be a $\lambda$-zone.
Note that $U(\Gamma) \cap \partial\Delta$ is a cycle, and the drawing $\Gamma' = \Gamma \cap (\Sigma - \Delta)$ is $2$-cell in $\Sigma$.
We say that $\Gamma'$ is the {\it drawing obtained from $\Gamma$ by clearing $\Delta$}.
We say that $\T'$ is a {\it tangle of order $\theta-4\lambda-2$ obtained by clearing $\Delta$} if $\T'$ is a tangle in $\Gamma'$ of order $\theta-4\lambda-2$, and
\begin{itemize}
	\item $\T'$ is respectful with a metric $m_{\T'}$, and 
	\item $\T'$ is conformal with $\T$, and
	\item if $x,y$ are atoms of $\Gamma$ and $x',y'$ are atoms of $\Gamma'$ with $x \subseteq x'$ and $y \subseteq y'$, then $m_\T(x,y) \geq m_{\T'}(x',y') \geq m_\T(x,y)-4\lambda-2$.
\end{itemize}

\begin{theorem}[{\cite[Theorem~(7.10)]{rs XII}}] \label{clean a zone}
Let $\Delta$ be a $\lambda$-zone.
If $\theta \geq 4\lambda+3$, then there exists a unique respectful tangle of order $\theta-4\lambda-2$ obtained by clearing $\Delta$.
\end{theorem}

\begin{theorem}[{\cite[Theorem~(9.2)]{rs XIV}}] \label{big zone contains ball}
Let $\Gamma$ be a $2$-cell drawing in a surface $\Sigma$, and let $\T$ be a respectful tangle in $\Gamma$ of order $\theta$.
Let $x$ be an atom of $\Gamma$, and $\lambda$ an integer with $2 \leq \lambda \leq \theta-4$.
Then there exists a $(\lambda+3)$-zone $\Delta$ around $x$ such that $x' \subseteq \Delta$ for every atom $x'$ of $\Gamma$ with $m_\T(x,x') \leq \lambda$.
\end{theorem}

\begin{lemma} \label{disjoint boundry of zone}
Let $\Gamma$ be a $2$-cell drawing in a surface, $z$ an atom, and $\T$ a respectful tangle in $\Gamma$ of order $\theta$.
Let $\lambda$ be a nonnegative integer, and let $C$ be the cycle of the boundary of a $\lambda$-zone around $z$.
If $\theta \geq \lambda+8$, then there exists a $(\lambda+7)$-zone $\Lambda$ around $z$ such that the cycle bounding $\Lambda$ is disjoint from $C$, and $\Lambda$ contains the $\lambda$-zone bounded by $C$.
\end{lemma}

\begin{pf}
For every atom $x$ of $\Gamma$, let $\Lambda_x$ be a $4$-zone around $x$ containing all atoms $y$ with $m_\T(x,y)\leq 1$, and let $\Delta_x$ be the closure of $\Lambda_x$, and let $C_x$ be the boundary cycle of $\Delta_x$.
For every $v \in V(C)$, since every region incident with $v$ has distance $1$ from $v$, $v$ is an interior point of $\Delta_v$.
Let $\Delta = \Delta' \cup \bigcup_{v \in V(C)} \Delta_v$, where $\Delta'$ is the open disk with the boundary $C$.
So $V(C)$ are interior points of $\Delta$.
By the triangle-inequality, for every $v \in V(C)$ and for every vertex $u$ in $\Delta_v$, $m_\T(z,u) \leq \lambda +4$.
Therefore, there exists a $(\lambda+7)$-zone $\Lambda$ around $z$ that contains $\Delta$ by Theorem \ref{big zone contains ball}.
Since any vertex in $C$ is an interior point of $\Delta$, it is an interior point of $\Lambda$, so $C$ is disjoint from the cycle that bounds $\Lambda$.
\end{pf}

\bigskip

Let $\Gamma$ be a 2-cell drawing in a surface having a respectful tangle $\T$.
Let $\Lambda$ be a $\lambda$-zone (with respect to $m_\T$) around some atom of $\Gamma$ for some nonnegative integer $\lambda$.
For every $v \in V(\Gamma) \cap \partial\Lambda$, a {\it loose component with respect to $(\Lambda,v,\T)$} is a component $L$ of $\Gamma-v$ such that some vertex of $L$ is adjacent to $v$, and there exists no separation $(A,B) \in \T$ with $V(A \cap B)=\{v\}$ and $V(B)=V(L) \cup \{v\}$; we call $v$ the {\it attachment} of $L$.
{\it A loose component with respect to $(\Lambda,\T)$} is a loose component with respect to $(\Lambda,v,\T)$ for some $v \in V(\Gamma) \cap \partial\Lambda$.

\begin{lemma} \label{making 2-cell}
Let $\Gamma$ be a $2$-cell drawing in a surface, $z$ an atom, and $\T$ a respectful tangle in $\Gamma$ of order $\theta$.
Let $\lambda$ be a nonnegative integer, and let $\Lambda$ be a $\lambda$-zone around $z$.
If $\theta \geq 4\lambda+31$, then there exists a $(\lambda+7)$-zone $\Lambda'$ around $z$ containing $\Lambda$ such that the following hold.
	\begin{enumerate}
		\item $\partial\Lambda'$ is obtained from $\partial\Lambda$ by adding pairwise disjoint paths in $\Gamma$, where each of them has both ends in $\partial\Lambda$ and for each of its internal vertices $v$, every edge incident in $\Gamma$ with $v$ is either contained in $\overline{\Lambda'}$ or incident with a vertex in a loose component with respect to $(\Lambda',v,\T)$.
		\item The drawing $\Gamma'$ obtained from $\Gamma$ by clearing $\Lambda'$, deleting $E(\Gamma) \cap \partial\Lambda'$ and all loose components with respect to $(\Lambda',\T)$ and deleting all resulting isolated vertices is 2-cell and has a respectful tangle $\T'$ of order at least $\theta-(4\lambda+30)$ conformal with $\T$ such that $m_{\T}(x',y') \geq m_{\T'}(x,y) \geq m_{\T}(x',y')-(4\lambda+30)$ for every two atoms $x,y$ of $\Gamma'$, where $x',y'$ are atoms of $\Gamma$ with $x' \subseteq x$ and $y' \subseteq y$.
	\end{enumerate}
\end{lemma}

\begin{pf}
By Lemma \ref{disjoint boundry of zone}, there exists a $(\lambda+7)$-zone $\Lambda_0$ such that $\Lambda_0 \supseteq \Lambda$ and $\partial\Lambda_0 \cap \partial\Lambda=\emptyset$.
Let $\Gamma_0$ be the drawing obtained from $\Gamma$ by clearing $\Lambda_0$.
Then $\Gamma_0$ has a respectful tangle $\T_0$ of order at least $\theta-(4\lambda+30)$ obtained from $\T$ by clearing $\Lambda_0$ conformal with $\T$, and $\Gamma_0$ is 2-cell.

Let $\Gamma_1$ be the drawing obtained from $\Gamma$ by clearing $\Lambda$ and deleting $E(\Gamma) \cap \partial\Lambda$ and isolated vertices.
Since $\Gamma_0$ is a subgraph of $\Gamma_1$, $\Gamma_1$ contains a respectful tangle $\T_1$ of order at least $\theta-(4\lambda+30)$ conformal with $\T$.
We may assume that $\Gamma_1$ is not 2-cell; otherwise we are done by choosing $\Lambda'=\Lambda$.

Since $\theta-(4\lambda+30) \geq 1$, there exists a component of $\Gamma_1$, say $Q_1$, such that $(\Gamma_1-V(Q_1),Q_1) \in \T_1$.
Since $\Gamma_0$ is 2-cell, $\Gamma_0$ is connected.
So $V(\Gamma_0) \subseteq V(Q_1)$.
Hence, every component of $\Gamma_1-V(Q_1)$ is contained in $\Lambda_0$.

Let $Q$ be a component of $\Gamma_1-V(Q_1)$ with $\lvert V(Q) \cap \partial\Lambda \rvert \geq 2$.
Hence there exists a path $P_Q$ in $Q$ on at least two vertices with ends contained in $\partial\Lambda$.
We choose $P_Q$ such that the closed set, denoted by $\Delta_Q$, bounded by $\partial\Lambda \cup P_Q$ with interior disjoint from $V(Q_1)$ is maximal.

Define $\Lambda'$ to be the open disk whose closure is $\overline{\Lambda} \cup \bigcup_{Q'} \Delta_{Q'}$, where the union is over all components $Q'$ of $\Gamma_1-V(Q_1)$ with $\lvert V(Q') \cap \partial\Lambda \rvert \geq 2$.
Hence $\partial\Lambda'$ is obtained from $\partial\Lambda$ by adding $\bigcup_{Q'}P_{Q'}$, where the union is over some components $Q'$ of $\Gamma_1-V(Q_1)$.
Clearly, $\Lambda'$ is bounded by a cycle $C'$ in $\Gamma$.
Note that those $P_{Q'}$ are pairwise disjoint, and for every internal vertex $v$ of some $P_{Q'}$, every edge incident in $\Gamma$ with $v$ is either contained in $\overline{\Lambda'}$ or incident with a vertex in a loose component with respect to $(\Lambda',v,\T)$.
Furthermore, the drawing $\Gamma'$ obtained from $\Gamma$ by clearing $\Lambda'$ and deleting $E(C')$, all loose components with respect to $(\Lambda',\T)$ and all isolated vertices is $Q_1$.
Since $\partial \Lambda \cap \partial \Lambda_0=\emptyset$, every region of $Q_1$ is either a region of $\Gamma$ or a subset of $\Lambda_0$.
Since $Q_1$ is connected, $Q_1$ and hence $\Gamma'$ is 2-cell.
Since $\Gamma_0$ is a subgraph of $\Gamma'$, there exists a respectful tangle $\T'$ in $\Gamma'$ of order at least $\theta-(4\lambda+30)$ conformal with $\T$, and $m_\T(x',y') \geq m_{\T'}(x,y) \geq m_{\T_0}(x'',y'') \geq m_{\T}(x',y')-(4\lambda+30)$ for all atoms $x,y$ of $\Gamma'$, where $x',y'$ are atoms of $\Gamma$ and $x'',y''$ are atoms of $\Gamma_0$ with $x' \subseteq x \subseteq x''$ and $y' \subseteq y \subseteq y''$.
Since $\Lambda' \subseteq \Lambda_0$, $\Lambda'$ is a $(\lambda+7)$-zone.
This proves the lemma.
\end{pf}

\bigskip

Let $\Sigma$ be a connected surface, and let $\Delta_1, ..., \Delta_t$ be pairwise disjoint closed disks in $\Sigma$.
Let $\Gamma$ be a drawing in $\Sigma$ such that $U(\Gamma) \cap \Delta_i = V(\Gamma) \cap \partial\Delta_i$ for $1 \leq i \leq t$.
Let $Z = \bigcup_{i=1}^t V(\Gamma) \cap \partial\Delta_i$.
We say that a partition $(Z_1, Z_2, ..., Z_p)$ of $Z$ satisfies the {\it topological feasibility condition} if there exist pairwise disjoint disks $D_1, D_2, ..., D_p$ in $\Sigma$ such that $D_j \cap (\bigcup_{i=1}^t \Delta_i) = Z_j$ for $1 \leq j \leq p$.

\begin{theorem}[{\cite[Theorem~(3.2)]{rs XII}}] \label{linkage on surface}
For every connected surface $\Sigma$ and all integers $t \geq 0$ and $z \geq 0$, there exists a positive integer $\theta \geq 1$ such that the following is true.
Let $\Delta_1, ..., \Delta_t$ be pairwise disjoint closed disks in $\Sigma$, and let $\Gamma$ be a $2$-cell drawing in $\Sigma$ such that $U(\Gamma) \cap \Delta_i = V(\Gamma) \cap \partial\Delta_i$ for $1 \leq i \leq t$.
Let $\lvert Z \rvert \leq z$, where $Z = \bigcup_{i=1}^t (V(\Gamma) \cap \partial\Delta_i)$, and let $(Z_1, Z_2,...,Z_p)$ be a partition of $Z$ satisfying the topological feasibility condition.
Let $\T$ be a respectful tangle of order at least $\theta$ in $\Gamma$ with metric $m_\T$ such that $m_\T(r_i,r_j) \geq \theta$ for $1 \leq i< j \leq t$, where $r_i$ is the region of $\Gamma$ meeting $\Delta_i$ for $1 \leq i \leq t$, and $V(\Gamma) \cap \partial\Delta_i$ is free for $1 \leq i \leq t$.
Then there are mutually disjoint connected subdrawings $\Gamma_1, \Gamma_2, ..., \Gamma_p$ of $\Gamma$ with $V(\Gamma_j) \cap Z = Z_j$ for $1 \leq j \leq p$.
\end{theorem}

\section{Excluding a subdivision of a fixed graph}
\label{sec:structures}

Let $G$ be a graph and $\T$ a tangle in $G$.
Given an integer $k$, a vertex $v$ of $G$ is said to be {\it $k$-free} (with respect to $\T$) if there is no $(A,B) \in \T$ of order less than $k$ such that $v \in V(A)-V(B)$.

Let $\Se=\Se_1 \cup \Se_2$ be a segregation of $\Se$ with $\Se_1 \cap \Se_2=\emptyset$ such that $\lvert \overline{\Omega} \rvert \leq 3$ for every $(S,\Omega) \in \Se_1$.
The {\it skeleton} of a proper arrangement $\alpha$ of $\Se$ in $\Sigma$ (with respect to $(\Se_1,\Se_2)$) is the drawing $\Gamma = (U,V)$ in $\Sigma$ with $V(\Gamma)=\{\alpha(v): v \in V(\Se)\}$ such that $U(\Gamma)$ consists of the boundary of $\alpha(S,\Omega)$ for each $(S,\Omega) \in \Se_1$ with $\lvert \bar{\Omega} \rvert = 3$, and a line in the boundary $\alpha(S',\Omega')$ with ends $\overline{\Omega'}$ for each $(S',\Omega') \in \Se_1$ with $\lvert \overline{\Omega'} \rvert=2$.
Note that we do not add any edges into the skeleton for $(S,\Omega)$ with $\lvert \bar{\Omega} \rvert \leq 1$ or $(S,\Omega) \in \Se_2$.
Furthermore, the skeleton of $\alpha$ is unique up to the choice of the line in $\alpha(S',\Omega')$ for each $(S',\Omega') \in \Se$ with $\lvert \overline{\Omega'} \rvert=2$, and the choices of those lines do not affect whether the skeleton is a 2-cell drawing or not.

\begin{lemma} \label{extend a vortex}
Let $t,\rho,\theta$ be nonnegative integers.
Let $G$ be a graph.
Let $\Se=\Se_1 \cup \Se_2$ with $\Se_1 \cap \Se_2=\emptyset$ be a segregation of $G$ such that $\lvert \overline{\Omega} \rvert \leq 3$ for every $(S,\Omega) \in \Se_1$.
Let $\alpha$ be a proper arrangement of $\Se$ with respect to $(\Se_1,\Se_2)$ of $G$ in a surface $\Sigma$.
Let $(S,\Omega) \in \Se$ be a $\rho$-vortex.
Let $G'$ be the skeleton of $\alpha$.
Let $\T'$ be a respectful tangle in $G'$ of order $\theta$.
If $G'$ is $2$-cell and $\theta \geq 4t+59$, then there exists a cycle $C$ such that the following hold.
	\begin{enumerate}
		\item $C$ bounds a $(t+14)$-zone $\Lambda$ in $G'$ around some vertex in $\bar{\Omega}$.
		\item $\Lambda$ contains every atom $x$ of $G'$ with $m_{\T'}(x,y) \leq t$ for some $y \in \bar{\Omega}$.
		\item The closure of $\Lambda$ contains $\alpha(S,\Omega)$.
		\item Let $S'$ be the union of $S''$ over all societies $(S'',\Omega'') \in \Se$ with 
			\begin{itemize}
				\item either $\alpha(S'',\Omega'') \subseteq \overline{\Lambda}$, or 
				\item $\lvert \overline{\Omega''} \rvert=2$ and $\alpha(S'',\Omega'') \cap E(C) \neq \emptyset$, or
				\item $\overline{\Omega''}$ contained in the union of some loose component with respect to $(\Lambda,\T')$ and its attachment, or 
				\item $\lvert \overline{\Omega''} \rvert=1$ and $\overline{\Omega''} \subseteq V(C)$.
			\end{itemize}
			Let $\overline{\Omega'}=V(C)-\{x \in V(C):$ every edge of $G'$ incident with $x$ is either contained in $\overline{\Lambda}$ or incident in $G'$ with a vertex in a loose component with respect to $(\Lambda,\T')\}$, and let $\Omega'$ be a cyclic ordering consistent with the cyclic ordering of $C$.
If every $(S'',\Omega'') \in \Se_2$ with $\alpha(S'',\Omega'') \subseteq \overline{\Lambda}$ is a $\rho_{S''}$-vortex for some nonnegative integer $\rho_{S''}$, then $(S',\Omega')$ is a $(\rho+t+8+ \sum_{S''} \rho_{S''})$-vortex, where the sum is over all societies $(S'',\Omega'') \in \Se_2-\{(S,\Omega)\}$ with $\alpha(S'',\Omega'') \subseteq \overline{\Lambda}$.
		\item Let $\Se_1^*=\Se_1-\{(S'',\Omega'') \in \Se_1: S'' \subseteq S'\}$ and $\Se^*_2=(\Se_2-\{(S'',\Omega'') \in \Se_2: S'' \subseteq S'\}) \cup \{(S',\Omega')\}$.
			If $m_{\T'}(x,y) \geq 3$ for every atom $x \subseteq \partial\Lambda$ and $y \in V(\overline{\Omega''})$ with $(S'',\Omega'') \in \Se_2^*-\{(S',\Omega')\}$, then $\Se^*$ is a segregation, and there exists a proper arrangement $\alpha^*$ of $\Se^*_1 \cup \Se^*_2$ with respect to $(\Se_1^*,\Se_2^*)$ such that the skeleton $G^*$ of $\alpha^*$ 
			\begin{itemize}
				\item is obtained from $G'$ by clearing $\Lambda$ and deleting some edges in $E(C)$ and all loose components with respect to $(\Lambda,\T')$ and deleting all resulting isolated vertices, 
				\item is 2-cell, and 
				\item has a respectful tangle $\T^*$ conformal with $\T'$ of order at least $\theta-4t-58$ such that $m_{\T'}(x',y') \geq m_{\T^*}(x,y) \geq m_{\T'}(x',y')-4t-58$ for all atoms $x,y$ of $G^*$, where $x',y'$ are atoms of $G'$ with $x' \subseteq x$ and $y' \subseteq y$.
			\end{itemize}
	\end{enumerate}
\end{lemma}

\begin{pf}
Let $y$ be a vertex in $\bar{\Omega}$.
By Theorem \ref{big zone contains ball}, there exists a $(t+5)$-zone $\Lambda'$ around $y$ in $G'$ such that $x \in \Lambda'$ for every atom $x$ of $G'$ with $m_{\T'}(x,y) \leq t+2$.
Since $m_{\T'}(y',y'') \leq 2$ for every two vertices $y',y''$ in $\bar{\Omega}$, $x \in \Lambda'$ for every atom $x$ of $G'$ with $m_{\T'}(x,z) \leq t$ for some $z \in \bar{\Omega}$.
Let $H$ be the drawing obtained from $G'$ by deleting every atom $x \in V(G')$ with $m_{\T'}(x,y) \leq t+2$.
It follows from \cite[Theorem~(8.10)]{rs XI} that
 $H$ has a region $f$ homeomorphic to an open disk that contains $\alpha(S,\Omega)$ and all deleted vertices.

In the rest of the proof, we fix a radial drawing of $G'$.

\noindent{\bf Claim 1:} For every vertex $v$ of $H$ incident with $f$, there exists a closed walk $\ell_v$ of length at most $2t+8$ in the radial drawing of $G'$ with $\{v,y\} \subseteq \ins(\ell_v) \subseteq \bar{f} \subseteq \overline{\Lambda'}$ such that $v$ is adjacent to only one vertex in $\ell_v$ and $V(\ell_v) \cap V(H) = \{v\}$.

\noindent{\bf Proof of Claim 1:}
Since $v$ is incident with $f$, there exists a path $P \subseteq \bar{f}$ of length two in the radial drawing of $G'$ containing $v$ and a vertex $v'$ of $G'-V(H)$ internally disjoint from $V(H)$.
As $m_{\T'}(v',y) \leq t+2$, there exists a closed walk $W_{v'}$ of length at most $2t+4$ in the radial drawing of $G'$ such that $\{v',y\} \subseteq \ins(W_{v'})$.
Note that $v \in V(H)$, so $m_{\T'}(v,y) >t+2$ and $\{v\} \not \subseteq \ins(W_{v'})$.
Hence, there exists a closed walk $\ell_v$ of length at most $2t+8$ in $W_{v'} \cup P$ with $\{v,y\} \subseteq \ins(\ell_v)$ and such that $v$ is adjacent to only one vertex in $\ell_v$.
Note that $m_{\T'}(u,y) \leq t+2$ for every $u \in V(W_{v'}) \cap V(G)$, so $\ins(W_{v'}) \subseteq \bar{f}$.
Hence $\ins(\ell_v) \subseteq \bar{f}$.
Since $\Lambda'$ contains all atoms $u$ of $G'$ with $m_{\T'}(u,y) \leq t+2$, $\bar{f} \subseteq \overline{\Lambda'}$.
$\Box$

By Claim 1, $m_{\T'}(y,v) \leq t+4$ for every vertex $v$ incident with $f$.
So there exists a $(t+7)$-zone $\Lambda''$ containing all vertices incident with $f$.
Hence $f \subseteq \Lambda''$.
By Lemma \ref{something is far}, there exists an edge $y'$ of $G'$ such that $m_{\T'}(y,y')=\theta>t+8$.
So $y' \not \in \Lambda''$.
Each cycle $Y$ contained in the boundary of $f$ bounds an open disk $\Delta_Y$ disjoint from $y'$.
Let $C_0$ be the cycle contained in the boundary of $f$ such that $\Delta_{C_0}$ is maximal among all cycles $Y$ contained in the boundary of $f$.
Then $f \subseteq \Delta_{C_0}$ and hence $\Delta_{C_0}$ contains $\alpha(S,\Omega)$ and all atoms $x$ with $m_{\T'}(x,y) \leq t+2$.
Since $\Delta_{C_0} \subseteq \Lambda''$, $C_0$ bounds a $(t+7)$-zone.

Let $S'_0$ be the union of $S''$ over all societies $(S'', \Omega'') \in \Se$ with $\alpha(S'',\Omega'')$ contained in the closure of the disk bounded by $C_0$ disjoint from $y'$.
We call such a society $(S'',\Omega'')$ an {\it inner vortex}.
Let $\Omega'_0$ be the cyclic ordering of $C_0$.
Since $(S,\Omega)$ is a $\rho$-vortex, for every two intervals $I,J$ that partition $\bar{\Omega}$, there exists $X_{I,J} \subseteq V(S)$ with $\lvert X_{I,J} \rvert \leq \rho$ such that there exists no path in $S-X_{I,J}$ from $I-X_{I,J}$ to $J-X_{I,J}$.
Similarly, since $(S'',\Omega'')$ is a $\rho_{S''}$-vortex for every $(S'',\Omega'') \in \Se_2-\{(S,\Omega)\}$, for every two intervals $I,J$ that partition $\overline{\Omega''}$, there exists $X_{I,J}^{S''} \subseteq V(S)$ with $\lvert X_{I,J}^{S''} \rvert \leq \rho_{S''}$ such that there exists no path in $S''-X_{I,J}^{S''}$ from $I-X_{I,J}^{S''}$ to $J-X_{I,J}^{S''}$.

\noindent{\bf Claim 2:} $(S'_0,\Omega'_0)$ is a $(\rho+t+6+\sum \rho_{S''})$-vortex, where the sum is over all inner vortices.

\noindent{\bf Proof of Claim 2:}
Let $I',J'$ be two intervals that partition $\overline{\Omega'}$, let $u,v$ be the first vertex in $I',J'$, respectively, under the ordering $\Omega'$, and let $\ell_u^*$ and $\ell_v^*$ be a closed walk $\ell_u$ and a closed walk $\ell_v$ mentioned in Claim 1, respectively.

Note that for each path $Q$ in $\ell^*_u \cup \ell^*_v$ and for each vertex $w \in V(Q)$ corresponding to the region containing $\alpha(S'',\Omega'')$ for some inner vortex $(S'',\Omega'')$ or vortex $(S'',\Omega'')=(S,\Omega)$, the edges in $Q$ incident with $w$ define a partition of $\overline{\Omega''}$ into two intervals $I'',J''$, and we denote $X^{S''}_{I'',J''}$ by $X^{S''}_{w,Q}$.

Assume that there exists a path $Q$ in $\ell^*_u \cup \ell^*_v$ from $u$ to $v$ on at most $2t+11$ vertices.
Let $X''=(V(Q) \cap V(G')) \cup \bigcup_{w \in V(Q)-V(G')}\bigcup_{S''}X_{w,Q}^{S''}$, where the last union is over all $(S'',\Omega'')$ such that either $(S'',\Omega'')=(S,\Omega)$ or $(S'',\Omega'')$ is an inner vortex with $\alpha(S'',\Omega'')$ corresponding to $w$.
Then there exists no path in $S'_0-X''$ from $I'-X''$ to $J'-X''$.
Note that $\lvert X'' \rvert \leq \lceil \frac{2t+11}{2}\rceil+\rho+\sum_{S''}\rho_{S''}$, where the sum is over all inner vortices $(S'',\Omega'')$.
Therefore, $(S'_0,\Omega'_0)$ is a $(\rho+t+6+\sum_{S''}\rho_{S''})$-vortex, where the sum is over all inner vortices $(S'',\Omega'')$.

So we may assume that there does not exist a path $Q$ in $\ell^*_u \cup \ell^*_v$ from $u$ to $v$ on at most $2t+11$ vertices.
In particular, $\ell^*_u$ is disjoint from $\ell^*_v$.
So one of $\ell^*_u$ and $\ell^*_v$ does not contain $y$.
By symmetry, we may assume that $\ell^*_u$ does not contain $y$.
Since $\{y\} \subseteq \ins(\ell^*_u)$, $\ell^*_u$ contains a cycle $C_u$ such that $C_u$ bounds an open disk $\Delta_u \subseteq \ins(\ell^*_u)$ containing $\{y\}$.
Since $\overline{\Delta_u} \subseteq \bar{f} \subseteq \overline{\Lambda'}$, $\{v\}$ is not contained in the closure of $\Delta_u$, for if it did, then, since $v$ belongs to the boundary of $f$, it would belong to the boundary of $\Delta_u$, contrary to the fact that $\ell^*_u$ and $\ell^*_v$ are disjoint. 
Since $\ell^*_u$ is disjoint from $\ell^*_v$, $\ell^*_v$ is disjoint from the closure of $\Delta_u$.
In particular, $y \not \in V(\ell^*_v)$.
Hence $\ell^*_v$ contains a cycle $C_v$ such that $C_v$ bounds an open disk $\Delta_v \subseteq \ins(\ell^*_v)$ containing $\{y\} \cup \Delta_u$.
Since $\overline{\Delta_v} \subseteq \bar{f} \subseteq \overline{\Lambda'}$, $\{u\}$ is not contained in the closure of $\Delta_v$, for if it did, then, since $u$ belongs to the boundary of $f$, it would belong to the boundary of $\Delta_v$, contrary to the fact that $\ell^*_u$ and $\ell^*_v$ are disjoint.
Since $\ell^*_u$ intersects the closure of $\Delta_v$ and the complement of the closure of $\Delta_v$, $\ell^*_u$ intersects $C_v \subseteq \ell^*_v$, a contradiction.
This proves the claim.
$\Box$

By Lemma \ref{making 2-cell}, there exists a $(t+14)$-zone $\Lambda'$ containing the $(t+7)$-zone bounded by $C_0$ satisfying the conclusion of Lemma \ref{making 2-cell}.
Define $C$ to be the cycle bounding $\Lambda'$.
So Conclusions 1-3 hold. 

Let $(S',\Omega')$ be the society mentioned in Conclusion 4.
To prove Conclusion 4, it suffices to prove the following claim.

\noindent{\bf Claim 3:} If $(S_0',\Omega'_0)$ is a $p$-vortex for some nonnegative integer $p$, then $(S',\Omega')$ is a $(p+2)$-vortex.

\noindent{\bf Proof of Claim 3:}
Suppose that $(S',\Omega')$ is not a $(p+2)$-vortex.
So there exist a partition $I,J$ of $\Omega'$ into cyclic intervals and $p+3$ disjoint paths $P_1,...,P_{p+3}$ in $S'$ from $I$ to $J$.
By the definition of $\Omega'$, we know $\overline{\Omega'} \subseteq \overline{\Omega'_0}$.
Let $I'$ be the minimal cyclic interval in $\Omega'_0$ containing $I$, and let $J'$ be the cyclic interval such that $I' \cup J'=\Omega_0'$.
Note that the ends of $I$ and $I'$ are the same.
By changing the indices, we may assume that $P_1,P_2,...,P_{p+1}$ do not intersect the ends of $I'$.
So for each $i$ with $1 \leq i \leq p+1$, some subpath $Q_i$ of $P_i$ in $S_0'$ is from $I'$ to $J'$.
This contradicts that $(S_0',\Omega_0')$ is a $p$-vortex and proves the claim.
$\Box$

Then Conclusion 4 follows from Claims 2 and 3.

Define $\Se_1^*$ and $\Se_2^*$ as mentioned in Conclusion 5.
Now we assume that $m_{\T'}(x,y) \geq 3$ for every $x \subseteq \partial\Lambda$ and $y \in \overline{\Omega''}$ with $(S'',\Omega'') \in \Se_2^*-\{(S',\Omega')\}$.  
Hence every loose component with respect to $(\Lambda,\T')$ does not intersect $\overline{\Omega''}$ for every $(S'',\Omega'') \in \Se_2^*-\{(S',\Omega')\}$.
So $\Se^*=\Se^*_1 \cup \Se^*_2$ is a segregation of $G$.
Then it is clear that there exists a proper arrangement $\alpha^*$ of $\Se^*$ with respect to $(\Se_1^*,\Se^*_2)$ such that the skeleton $G^*$ of $\alpha^*$ can be obtained from $G'$ by deleting some edges in $E(C)$ and all loose components with respect to $(\Lambda,\T')$, and deleting all resulting isolated vertices.
Conclusion 2 of Lemma \ref{making 2-cell} implies that $G^*$ is 2-cell.
Furthermore, Conclusion 2 of Lemma \ref{making 2-cell} implies that there exists a respectful tangle $\T^*$ in $G^*$ conformal with $\T'$ of order at least $\theta-(4(t+7)+30)=\theta-4t-58$ such that $m_{\T'}(x',y') \geq m_{\T^*}(x,y) \geq m_{\T'}(x',y')-(4t+58)$ for all atoms $x,y$ of $G^*$, where $x',y'$ are atoms of $G'$ with $x' \subseteq x$ and $y' \subseteq y$.
This proves Conclusion 5.
\end{pf}

\begin{lemma} \label{free in vortex}
Let $d \geq 3$, and let $\kappa,h, h_1,h_2,...,h_\kappa, \rho, \theta''$ be nonnegative integers.
Then there exist integers $\theta_0(d,h,\rho,\kappa,\theta'')$, $\beta(d,h,\rho)$ and $f(d,h,\rho,\kappa)$ such that the following holds.
Suppose that
	\begin{enumerate}
		\item $G$ is a graph and $\T$ is a tangle in $G$, and 
		\item $\tau$ is a proper arrangement of a $\T$-central segregation $\Se$ of $G$ with respect to $(\Se_1,\Se_2)$ in a surface $\Sigma$, and 
		\item $G'$ is the skeleton of $\tau$, $G'$ is $2$-cell and is a minor of $G$, and $\T'$ is a respectful tangle in $G'$ of order $\theta$, for some $\theta \geq \theta_0$, such that $\T'$ is conformal with $\T$, and 
		\item let $X \subseteq V(G)$ and let $(S_1, \Omega_1), ..., (S_\kappa, \Omega_\kappa)$ be societies in $\Se$, where each $(S_i,\Omega_i)$ is a $\rho$-vortex such that $V(S_i) \cap X$ contains at least one $d$-free vertex with respect to $\T$, such that for every $1 \leq i < j \leq \kappa$, and for every $x \in \overline{\Omega_i}$ and $y \in \overline{\Omega_j}$, $m_{\T'}(x,y) \geq 2f+1$, and
		\item $m_{\T'}(x,y) \geq f+1$, for every $x \in \overline{\Omega_i}$ with $1 \leq i \leq \kappa$, and for every $y \in \bar{\Omega}$ with $(S,\Omega) \in \Se-\{(S_i,\Omega_i)\}$ and $\lvert \bar{\Omega} \rvert >3$, and 
		\item $h_i \leq h$ for $1 \leq i \leq \kappa$.
	\end{enumerate}
Then there exist $Z_1, Z_2, ..., Z_\kappa \subseteq V(G)$, a subdrawing $G''$ of $G'$ and a tangle $\T''$ in $G''$ obtained from $G'$ and $\T'$ by clearing at most $\kappa$ $f$-zones $\Lambda_1,...,\Lambda_\kappa$ such that $\T''-\bigcup_{i=1}^\kappa Z_i$ has order at least $\theta''$ and is conformal with $\T'$, and for every $i \in \{1,2,...,\kappa\}$, either 
	\begin{enumerate}
		\item $h_i \geq 2$, $\Lambda_i = \emptyset$ and $\lvert Z_i \rvert \leq \beta$ such that every vertex in $(V(S_i)-Z_i) \cap X$ is not $d$-free with respect to $\T-\bigcup_{i=1}^\kappa Z_i$, or
		\item $Z_i = \emptyset$, and $\Lambda_i$ is an $f$-zone in $G'$ around a vertex in $\overline{\Omega_i}$ with the boundary cycle $Y_i$, and there exist $h_i$ subsets $A_{i,1}, A_{i,2}, ..., A_{i,h_i}$ of $Y_i$ such that the following hold.
		\begin{enumerate}
			\item $V(S_i) \subseteq \Lambda_i$.
			\item Each $A_{i,j}$ has size $d$ and $\bigcup_{j=1}^{h_i}A_{i,j}$ is free in $G''$ with respect to $\T''-\bigcup_{i=1}^\kappa Z_i$.
			\item For each $1 \leq j \leq h_i$, there exists a minimum interval $I_j$ of $Y_i$ containing $A_{i,j}$ such that $I_{j'} \cap I_{j''}=\emptyset$ for every $1 \leq j' < j'' \leq h_i$. 
			\item There exist $v_{i,1}, v_{i,2},...,v_{i,h_i} \in \Lambda_i \cap X$ such that there are $h_i$ disjoint $d$-spiders in $G$ contained in $\Lambda_i \cup \bigcup_{j=1}^{h_i}A_{i,j}$, where the $j$-th spider is from $v_{i,j}$ to $A_{i,j}$.
		\end{enumerate}
	\end{enumerate}
\end{lemma}

\begin{pf}
Define $k'$ to be the value $k'(h,\rho)$ mentioned in Theorem \ref{nest buffer makes good}, and let $\beta(d,h,\rho) = 2(k'd+1)^{d+1}$.
Define $s'=s'_{\ref{nest buffer makes good}}(h,d,4hd+2\kappa\beta+3,\rho)+2hd+\kappa\beta$, where $s'_{\ref{nest buffer makes good}}$ is the value $s'$ mentioned in Theorem \ref{nest buffer makes good}.
Let $f(d,h,\rho,\kappa)=19+10s'$ and $\theta_0(d,h,\rho,\kappa,\theta'') = \theta'' + \kappa(4f+\beta+2)$.
Let $i \in \{1,2,...,\kappa\}$ be fixed.
For simplicity, we denote $(S_i,\Omega_i)$ by $(S,\Omega)$, and let $v_S$ be a vertex in $\bar{\Omega}$.

We may assume that $\kappa \geq 1$ and $h \geq 1$, for otherwise this lemma is obvious.
In particular, $\theta_0 \geq 6$.
By Theorem \ref{big zone contains ball}, there exists a $5$-zone $\Lambda_{S,0}'$ in $G'$ around $v_S$ such that $\Lambda_{S,0}'$ contains all atoms $y$ of $G'$ with $m_{\T'}(v_S,y) \leq 2$.
Note that every vertex in $\bar{\Omega}$ has distance at most $2$ from $v_S$ with respect to the metric $m_{\T'}$, so $\Lambda_{S,0}'$ contains $\tau(S,\Omega)$.
Let $\Lambda_{S,0}$ be a $19$-zone in $G'$ around $v_S$ such that $\Lambda_{S,0}$ satisfies Lemma \ref{extend a vortex} and contains $\Lambda_{S,0}'$. 
Let $(G_{S,0}, \Omega_{S,0})$ be the society $(S',\Omega')$ mentioned in Lemma \ref{extend a vortex} by taking $\Lambda$ to be $\Lambda_{S,0}$.
Note that Lemma \ref{extend a vortex} ensures that $(G_{S,0},\Omega_{S,0})$ is a $(\rho+13)$-vortex.

For $1 \leq j \leq s'$, let $\Lambda_{S,j}$ be a $(19+10j)$-zone around $v_S$ such that $\Lambda_{S,j}$ contains every vertex $x$ of $G''$ with $m_{\T'}(x,v_S) \leq 19+10(j-1)$ and $\partial\Lambda_{S,j} \cap \partial\Lambda_{S,j-1} = \emptyset$.
Note that the existence of $\Lambda_{S,j}$ follows from Lemmas \ref{big zone contains ball} and \ref{disjoint boundry of zone}.
Let $C_{S,j}$ be the boundary cycle of $\Lambda_{S,j}$ for $1 \leq j \leq s'$.
Let $\Lambda_S=\Lambda_{S,s'}$.
Let $G_S$ be the union of $S'$ over all societies $(S',\Omega')$ with $\tau(S',\Omega')$ contained in the closure of $\Lambda_S$, and let $\Omega_S$ be the cyclic ordering on the boundary cycle of $\Lambda_S$.
So $(G_S, \Omega_S)$ is a composition of a $(\rho+13)$-vortex $(G_{S,0}, \Omega_{S,0})$ with a rural neighborhood which has a presentation with an $s'$-nest $(C_{S,1}, C_{S,2}, ..., C_{S,s'})$.

Let $h'_i = k'$ if $h_i \neq 1$, and $h'_i=1$ if $h_i=1$.
Let $X_S$ be the set of $d$-free vertices in $V(S) \cap X$ with respect to $\T$.
Note that $X_S \neq \emptyset$ by assumption.
By Theorem \ref{vertex spider set}, either there exist $h'_i$ disjoint $d$-spiders from $X_S$ to $\overline{\Omega_S}$, or there exists $W_S \subseteq V(G_S)$ with $\lvert W_S \rvert \leq 2(h_i'd+1)^{d+1} \leq \beta$ such that every $d$-spider from $X_S$ to $\overline{\Omega_S}$ intersects $W_S$.

We first assume that $h_i>1$ and the set $W_S$ mentioned above exists.
Then for every vertex $v \in X_S-W_S$, there exists a separation $(A,B)$ of $G_S-W_S$ of order less than $d$ such that $v \in V(A)-V(B)$ and $\overline{\Omega_S}-W_S \subseteq V(B)$.
So there exists a separation $(A',B')$ of $G-W_S$ with $V(A \cap B)=V(A' \cap B')$ such that $V(A')=V(A)$ and $V(B')=V(B) \cup (V(G)-V(G_S))$.
Note that $(A',B') \in \T-W_S$; otherwise since $\T'$ is conformal with $\T$, there exists a separation $(B'',A'') \in \T'-W_S$, where $V(A'')=V(A') \cap V(G')$ and $V(B'')=V(B') \cap V(G')$, but $A''$ is contained in a $(19+10s')$-zone around $v_S$ in $G'$, a contradiction.
So every vertex in $X_S-W_S$ is not $d$-free in $\T-W_S$.
Hence every vertex in $(V(S)-W_S) \cap X$ is not $d$-free in $\T-W_S$.
In other words, the first statement of this theorem holds by taking $\Lambda_i=\emptyset$ and $Z_i = W_S$.

When $h_i=1$, the former case holds by Menger's Theorem and the fact that $V(S) \cap X$ contains a $d$-free vertex.
Therefore, we may assume that there exist $h_i'$ disjoint $d$-spiders from $X_S$ to $\overline{\Omega_S}$. 

Define $Z_i$ to be the empty set.
Let $D_{i,1}, D_{i,2}, ..., D_{i,h'_i}$ be disjoint $d$-spiders from $X_{S_i}$ to $\overline{\Omega_{S_i}}$.
Apply Theorem \ref{nest buffer makes good} by taking $(S,\Omega)=(G_{S_i},\Omega_{S_i})$, $(S_0,\Omega_0)=(S_i,\Omega_i)$ and $D_j = D_{i,j}$ for $1 \leq j \leq h'_i$. 
There exist pairwise disjoint $d$-spiders $D'_{i,1}, D'_{i,2}, ..., D'_{i,h_i}$ from $X_{S_i}$ to $V(C_{S_i,s'})$, a $(4hd+2\kappa\beta+3)$-nest $(N_{S_i,1},...,N_{S_i,4hd+2\kappa\beta+3})$ and intervals $I_{i,1},I_{i,2},...,I_{i,h_i}$ of $C_{S_i,s'}$ satisfying the conclusions of Theorem \ref{nest buffer makes good}.
For every $1 \leq j \leq h_i$, since each $D'_{i,j}$ is perpendicular to $(N_{S_i,1},...,N_{S_i,4hd+2\kappa\beta+3})$, there exists a set $A_{i,j}$ of $h_id$ vertices in $D'_{i,j} \cap V(N_{S_i,1})$ such that there exist $h_id$ disjoint paths from $A_{i,j}$ to $V(C_{S_i,s'})$, but there exists no path from the head of $D'_{i,j}$ to $V(N_{S_i,1})$ in $D'_{i,j}-A_{i,j}$.
Note that $N_{S_i,1}$ is contained in the disk bounded by $C_{S_i,s'}$ which bounds an $f$-zone, so $N_{S_i,1}$ is the boundary of an $f$-zone.
Define $\Lambda_i$ to be the $f$-zone bounded by $N_{S_i,1}$.
Define $G''$ to be the drawing and $\T''$ to be the tangle obtained from $G'$ and $\T'$, respectively, by clearing $\bigcup_{i=1}^\kappa \Lambda_i$. 
Note that $\T''-\bigcup_{i=1}^\kappa Z_i$ has order $\theta - \kappa\beta - \kappa(4f+2) \geq \theta''$  and is conformal with $\T'-\bigcup_{i=1}^\kappa Z_i$.
On the other hand, by planarity, for every $1 \leq j \leq h_i$, there exists an interval $J_{i,j}$ of $N_{S_i,1}$ containing $A_{i,j}$, such that $J_{i,j} \cap J_{i,j'} = \emptyset$ for every $j' \neq j$.

To prove this lemma, it is sufficient to show that $\bigcup_{j=1}^{h_i}A_{i,j}$ is free with respect to $\T''-\bigcup_{j=1}^\kappa Z_j$.
Suppose that $\bigcup_{j=1}^{h_i}A_{i,j}$ is not free with respect to $\T''-\bigcup_{j=1}^\kappa Z_j$ for some $i$, then there exists $(A',B') \in \T''-\bigcup_{j=1}^\kappa Z_j$ with order less than $dh_i$ such that $\bigcup_{j=1}^{h_i}A_{i,j} \subseteq V(A')$.
Let $(A,B) \in \T''$ with $V(A)=V(A') \cup \bigcup_{j=1}^\kappa Z_j$ and $V(B)=V(B') \cup \bigcup_{j=1}^\kappa Z_j$.
We assume that $A$ is as small as possible, so $m_{\T''}(\bigcup_{j=1}^{h_i}A_{i,j},u) < dh+\kappa\beta$ for every $u \in V(A)$ by Theorem \ref{A distance}.

We claim that $\overline{\Omega_{S_i}} \subseteq V(B)-V(A)$.
Suppose to the contrary that there exists $u \in \overline{\Omega_{S_i}} \cap V(A)$.
So there exists a closed walk $W$ of a radial drawing of $G''$ with length less than $2dh_i+2\kappa\beta$ and $\{u,v\} \subseteq \ins(W)$ for some vertex $v \in \bigcup_{j=1}^{h_i}A_{i,j}$.
Since $\{v\} \subseteq \ins(W)$, if $W$ contains the vertex of a radial drawing of $G''$ corresponding to the region bounded by $N_{S_\ell,1}$ for any $\ell \neq i$, then $W \cap V(N_{s_\ell,j}) \neq \emptyset$ for every $2 \leq j \leq 4hd+2\kappa\beta+2$, so $\lvert W \rvert \geq 4hd+2\kappa\beta+1$, a contradiction.
So $W$ does not contain the vertex of a radial drawing of $G''$ corresponding to the region bounded by $N_{S_\ell,1}$ for any $\ell \neq i$.
Similarly, since $\{u\} \subseteq \ins(W)$, if $W$ contains the vertex of a radial drawing of $G''$ corresponding to the region bounded by $N_{S_i,1}$, then $W \cap V(N_{s_i,j}) \neq \emptyset$ for every $2 \leq j \leq 4hd+2\kappa\beta+2$, a contradiction.
So $W$ does not contain the vertex of a radial drawing of $G''$ corresponding to the region bounded by $N_{S_i,1}$.
Therefore, $W$ is a closed walk in a radial drawing of $G'$.
This implies that $m_{\T'}(u,v)<dh_i+\kappa\beta$.
So $u$ is a vertex in $C_{S_i,s'}$ contained in the disk bounded by $C_{S_i,dh+\kappa\beta}$ containing $C_{S_i,1}$, a contradiction.

Therefore, $\overline{\Omega_{S_i}} \subseteq V(B')-V(A')$.
However, there exist $dh_i$ disjoint paths in $G''$ from $\bigcup_{j=1}^{h_i}A_{i,j}$ to $\overline{\Omega_{S_i}}$, a contradiction.
So $\bigcup_{j=1}^{h_i}A_{i,j}$ is free with respect to $\T''-\bigcup_{j=1}^\kappa Z_j$ for every $i$.
This proves the lemma.
\end{pf}

\begin{lemma} \label{free should far apart}
Let $d \geq 3,h$ be positive integers, and let $\Sigma$ be a surface.
Then there exist integers $\theta(d,h, \Sigma), \phi(d,h,\Sigma)$ such that if $G$ is a $2$-cell drawing in $\Sigma$, $\T$ is a respectful tangle in $G$ of order at least $\theta$, and $G$ contains $h$ $d$-free vertices $v_1, v_2, ..., v_h$ with $m_\T(v_i, v_j) > \phi$ for $1 \leq i < j \leq h$, then $G$ admits an $H$-subdivision with branch vertices $v_1,v_2,...,v_h$, for every graph $H$ of order $h$ and of maximum degree at most $d$ embeddable in $\Sigma$.
\end{lemma}

\begin{pf}
Let $H$ be a graph of order $h$ and of maximum degree at most $d$ embeddable in $\Sigma$.
Let $\theta_{\ref{linkage on surface}}$ be the positive integer $\theta$ mentioned in Theorem \ref{linkage on surface} by taking $t=h$ and $z=dh$.
Note that $(\{v_i\},\{v_i\})$ is a $0$-vortex for every $i$.
For $1 \leq i \leq h$, let $\Lambda_i$ be the $12$-zone around $v_i$ of $G$ mentioned in Lemma \ref{extend a vortex} such that $\Lambda_i$ contains $v_i$ and all its neighbors, and let $S_i$ be the subgraph of $G$ contained in the closure of $\Lambda_i$, and let $\overline{\Omega_i} = \partial\Lambda_i \cap V(G)$ with the cyclic order defined by the boundary cycle of $\Lambda_i$.
So $(S_i,\Omega_i)$ is a $24$-vortex by Lemma \ref{extend a vortex}.
Let $\theta' = \theta_{\ref{free in vortex}}(d,1,24,h,\theta_{\ref{linkage on surface}})$, $\beta = \beta_{\ref{free in vortex}}(d,1,24)$ and $f=f_{\ref{free in vortex}}(d,1,24,\kappa)$, where $\theta_{\ref{free in vortex}}$, $\beta_{\ref{free in vortex}}$ and $f_{\ref{free in vortex}}$ are the numbers $\theta_0,\beta,f$ mentioned in Lemma \ref{free in vortex}.
Define $\theta = \theta'+h(4f+2)+2f+1$ and $\phi = \theta_{\ref{linkage on surface}}+h(4f+2)+2f+1$. 

Applying Lemma \ref{free in vortex} by taking $\kappa = h$, $h_i=1$ for $1 \leq i \leq h$, $\rho=24$, $\theta''=\theta_{\ref{linkage on surface}}$, and $\Se$ the segregation consisting of $(S_1,\Omega_1), (S_2,\Omega_2),...,(S_h,\Omega_h)$ and the societies in which each of them consists of exactly one edge that is not in $\bigcup_{i=1}^h S_i$, we obtain the desired subgraph $G''$ with a respectful tangle $\T''$, and $A_{i,1}$ for $1 \leq i \leq h$, such that every $A_{i,1}$ is free with respect to $\T''$, as mentioned in the conclusion of Lemma \ref{free in vortex}.
Then for every $x \in A_{i,1}$ and $y \in A_{j,1}$ for some $i \neq j$, we have that $m_{\T''}(x,y) \geq \theta_{\ref{linkage on surface}}$ by Theorem \ref{clean a zone}.

For $1 \leq i \leq h$, let $\Delta_i$ be a closed disk in $\Sigma$ contained in the closure of $\Lambda_i$ such that $\Delta_i \cap G'' = A_{i,1}$.
Since $H$ can be embedded in $\Sigma$, we can partition $\bigcup_{i=1}^h A_{i,1}$ and apply Theorem \ref{linkage on surface} to obtain a linear forest so that an $H$-subdivision in $G$ can be obtained by concatenating these linear forests and $h$ disjoint $d$-spiders $D_1,D_2,...,D_h$, where each $D_i$ is from $v_i$ to $A_{i,1}$ contained in $S_i$.
We obtain an $H$-subdivision in $G$ with branch vertices $v_1,v_2,...,v_h$.
\end{pf}

\begin{lemma} \label{make segregation central}
Let $\rho$ be an integer, $G$ a graph, $\T$ a tangle in $G$ of order at least $2\rho+2$, and $\Se$ a segregation of $G$.
If $(S,\Omega) \in \Se$ is a $\rho$-vortex and there exists no $(A,B) \in \T$ of order at most $2\rho+1$ such that $B \subseteq S$, then there exists no $(A',B') \in \T$ of order at most the half of the order of $\T$ such that $B' \subseteq S$.
\end{lemma}

\begin{pf}
Suppose that there exists $(A,B) \in \T$ of order at most the half of the order of $\T$ such that $B \subseteq S$.
Let $\bar{\Omega} = v_1,v_2,...,v_n$ in order, where $n=\lvert \bar{\Omega} \rvert$.
We may assume that every $v_i$ is adjacent to a vertex in $G-V(S)$, for otherwise we may remove it from $\bar{\Omega}$.
As $(S,\Omega)$ is a $\rho$-vortex, by Theorem 8.1 in \cite{rs IX}, there exists a path-decomposition $(P,\X)$ of $S$ of adhesion at most $\rho$ such that the $i$-th bag $X_i$ of $(P,\X)$ contains $v_i$ for every $1 \leq i \leq n$.
For every subgraph $H$ of $S$, we define $(A_H,B_H)$ to be the separation of $G$ with minimum order such that $A_H=H$.
In particular, for $1 \leq i \leq n$, $(A_{S[X_i]},B_{S[X_i]})$ has order at most $2\rho+1$, so $(A_{S[X_i]},B_{S[X_i]}) \in \T$.
For $1 \leq i \leq n$, define $(A_i,B_i) = (A \cup A_{S[\bigcup_{j=1}^i X_j]}, B \cap B_{S[\bigcup_{j=1}^i X_j]})$.
Note that if $v_j \in V(B)$ for some $j$, then $v_j \in V(A)$ since $B \subseteq S$ and $v_j$ is adjacent to a vertex in $G-V(S)$. 
So the order of $(A_i,B_i)$ is at most $\lvert V(A) \cap V(B) \rvert + \lvert V(A_{S[\bigcup_{j=1}^i X_j]}) \cap V(B_{S[\bigcup_{j=1}^i X_j]}) \cap (V(B)-V(A)) \rvert \leq \lvert V(A) \cap V(B) \rvert + \lvert (\bigcup_{j=1}^i X_j) \cap (\bigcup_{j=i+1}^nX_j) \rvert \leq \lvert V(A) \cap V(B) \rvert + \rho$.
Since the order of $(A,B)$ is at most the half of the order of $\T$, and the order of $\T$ is greater than $2\rho$, either $(A_i,B_i) \in \T$ or $(B_i,A_i) \in \T$ by the first tangle axiom.
Let $(A_0,B_0)=(A,B)$.
We shall prove that $(A_i,B_i) \in \T$ for $0 \leq i \leq n$ by induction on $i$.

When $i=0$, $(A_0,B_0)=(A,B) \in \T$.
Assume that $(A_i,B_i) \in \T$ for some $i$.
Suppose that $(B_{i+1},A_{i+1}) \in \T$.
But $(A_i,B_i),(A_{S[X_{i+1}]},B_{S[X_{i+1}]}) \in \T$, and $B_{i+1} \cup A_i \cup S[X_{i+1}] = G$, a contradiction.
This proves that $(A_i,B_i) \in \T$ for every $0 \leq i \leq n$.

Furthermore, $(A_n,B_n) = (A \cup S, B \cap B_S)$.
Recall that $V(B \cap B_S) \subseteq V(B) \cap \bar{\Omega} \subseteq V(A) \cap V(B)$, so $\lvert V(B_n) \rvert \leq \lvert V(A) \cap V(B) \lvert$.
Since $(B_n, G-E(B_n))$ has order less than the order of $\T$, $(B_n, G-E(B_n)) \in \T$ by the first and third tangle axioms.
However, $A_n \cup B_n=G$, contradicting the second tangle axiom.
This completes the proof.
\end{pf}

\begin{lemma} \label{make vortex away}
For a positive nondecreasing function $\phi$ and integers $\rho, \lambda, \kappa,k, \theta^*, d$ with $d \geq 4$, there exist integers $\theta, \rho^*$ such that the following is true.
Assume that $G$ is a graph, $X$ is a subset of $V(G)$, $\T$ is a tangle in $G$, and $\Se= \Se_1 \cup \Se_2$ is a $\T$-central segregation that has a proper arrangement $\tau$ in a surface $\Sigma$ such that the following hold.
	\begin{enumerate}
		\item For every $(S,\Omega) \in \Se_1$, $\lvert \bar{\Omega} \rvert \leq 3$, and for every $x \in \bar{\Omega}$, there exist $\lvert \bar{\Omega} \rvert-1$ paths in $S$ from $x$ to $\bar{\Omega}-\{x\}$ intersecting only in $\{x\}$.
		\item $\lvert \Se_2 \rvert \leq \kappa$.
		\item $(S,\Omega)$ is a $\rho$-vortex for every $(S,\Omega) \in \Se_2$.
		\item The skeleton $G'$ of $\Se$ is a minor of $G$, is $2$-cell embedded in $\Sigma$ and has a respectful tangle $\T'$ of order at least $\theta$ conformal with $\T$.
		\item There exist $k$ $\lambda$-zones $\Lambda_1,\Lambda_2,...,\Lambda_{k}$ in $G'$ around some vertices of $G'$ with respect to the metric $m_{\T'}$ such that every $d$-free vertex of $G'$ with respect to $\T'$ contained in $X$ is contained in $\bigcup_{i=1}^{k} \Lambda_{i}$.
	\end{enumerate}
Then there exist a $\T$-central segregation $\Se^*=\Se_1^* \cup \Se_2^*$ of $G$ and a proper arrangement $\tau^*$ of $\Se^*$ with respect to $(\Se^*_1,\Se^*_2)$ in $\Sigma$ such that the following hold.
	\begin{enumerate}
		\item $\Se_1^* \subseteq \Se_1$; in particular, $\lvert \bar{\Omega} \rvert \leq 3$ for every $(S,\Omega) \in \Se_1^*$.
		\item $\lvert \Se_2^* \rvert \leq \kappa+k$ and $\bigcup_{(S,\Omega) \in \Se_2} S \cup \bigcup_{(S',\Omega') \in \Se, \alpha(S',\Omega') \subseteq \bigcup_{i=1}^k \overline{\Lambda_i}} S' \subseteq \bigcup_{(S,\Omega) \in \Se_2^*} S$.
		\item There exists an integer $\rho'$ with $\rho' \leq \rho^*$ such that $(S,\Omega)$ is a $\rho'$-vortex for every $(S,\Omega) \in \Se_2^*$.
		\item The skeleton $G^*$ of $\tau^*$ is a minor of $G$, is $2$-cell embedded in $\Sigma$ and has a respectful tangle $\T^*$ of order at least $\theta^*+\phi(\rho^*)+2\rho^*$ conformal with $\T$.
		\item There is no $d$-free vertex of $G^*$ with respect to $m_{\T^*}$ contained in $X$.
		\item $m_{\T^*}(x,y) \geq \phi(\rho')$ for every atoms $x,y$ of $G^*$ with $x \in S_x, y \in S_y$ for different members $(S_x,\Omega_x), (S_y,\Omega_y) \in \Se_2^*$. 
	\end{enumerate}
\end{lemma}

\begin{pf}
Note that each society that consists of a single vertex is a $0$-vortex.
For each $i$, applying Lemma \ref{extend a vortex} by taking $(S,\Omega)=(\{v_i\},\{v_i\})$, where $v_i$ is a vertex of $G'$ such that $\Lambda_i$ is a $\lambda$-zone around $v_i$ (we will choose $\theta$ to be greater than $4\lambda+59$ so that Lemma \ref{extend a vortex} is applicable), we can find a $(\lambda+14)$-zone $\Lambda_i'$ containing $\Lambda_i$ such that the society $(S',\Omega')$ mentioned in Lemma \ref{extend a vortex} by taking $\Lambda=\Lambda_i$ is a $(\lambda+8+\kappa\rho)$-vortex. 
Therefore, we can replace $\Lambda_i$ by $\Lambda'_i$ so that we may assume that every $\Lambda_i$ is a $\lambda'$-zone and the subgraph of $G$ inside the disk $\Lambda_i$ is a $\lambda'$-vortex $(S,\Omega)$, where $\lambda'=\lambda+14+\kappa\rho$.
Similarly, for each $(S,\Omega) \in \Se_2$, there exists a $16$-zone $\Lambda_S$ containing the disk $\tau(S,\Omega)$, and the society $(S',\Omega')$ mentioned in Lemma \ref{extend a vortex} by taking $\Lambda=\Lambda_S$ is a $(\kappa\rho+10)$-vortex. 

Let $\C = \{\Lambda_i, \Lambda_S: 1 \leq i \leq k, (S,\Omega) \in \Se_2\}$, and let $\lambda_\C$ be the minimum $t$ such that every member of $\C$ is a $t$-zone.
For each member $\Lambda$ of $\C$, let $(S_\Lambda,\Omega_\Lambda)$ be the $(S',\Omega')$ mentioned in Lemma \ref{extend a vortex} by taking $\Lambda=\Lambda$. Let $M_\C$ be the minimum such that $(S_\Lambda,\Omega_\Lambda)$ are $M_\C$-vortices for all members $\Lambda$ of $\C$.
Let $\rho_0=M_\C$ and $\lambda_0=\lambda_\C$.
Note that $\lvert \C \rvert \leq k+\kappa$, $\rho_0 \leq \max\{\lambda', \kappa\rho+10\}$, and $\lambda_0 \leq \max\{\lambda', 16\}$.
For $i \geq 1$, let $t_i=\phi(\rho_{i-1})+(4\lambda_{i-1}+58)(k+\kappa)+2$, $\lambda_i=t_i+14$ and $\rho_i=(k+\kappa-i+1)\rho_{i-1}+t_i+8$.
Then we consecutively test whether there exist two atoms of $G'$ in different members of $\C$ with distance less than $\phi(M_\C)+(4\lambda_\C+58)\lvert \C \rvert$ under the metric $m_{\T'}$, and if such two nearby vortices exist, then we do the following.
Find the $(t_{k+\kappa+1-\lvert \C \rvert}+14)$-zone $\Lambda$ mentioned in the conclusion of Lemma \ref{extend a vortex} containing these two nearby members of $\C$, remove these two members from $\C$ and add $\Lambda$ into $\C$, and update $M_\C$ and $\lambda_\C$.
Since $\lvert \C \rvert$ decreases in each step, this process will terminate within $\kappa+k$ steps.
Furthermore, during the process, $\lambda_\C \leq \lambda_{k+\kappa-\lvert \C \rvert}$ and $M_\C \leq \rho_{k+\kappa-\lvert \C \rvert}$.
Therefore, when the process terminates, each member $\Lambda$ of $\C$ is a $\lambda_\C$-zone in $G'$ with $\lambda_\C \leq \lambda_{k+\kappa}$ and defines an $M_\C$-vortex $(S_\Lambda,\Omega_\Lambda)$ with $M_\C \leq \rho_{k+\kappa}$, and the distance between any two members of $\C$ is at least $\phi(M_\C)+(4\lambda_\C+58)\lvert \C \rvert$ under the metric $m_{\T'}$.
Let $\rho^*=\rho_{k+\kappa}$ and $\lambda^*=\lambda_{k+\kappa}$.
We define $\theta = \theta^*+\phi(\rho^*)+2\rho^*+1+(4\lambda^*+58)(\kappa+k)$. 

Define $\Se^*_1=\Se_1-\{(S'',\Omega'') \in \Se_1: S'' \subseteq \bigcup_{\Lambda \in \C}S_\Lambda\}$ and $\Se^*_2=\{(S_\Lambda,\Omega_\Lambda): \Lambda \in \C\}$.
Then Conclusions 1-3 hold.
Since atoms belonging to different members of $\Se_2^*$ have distance at least $3+(\lvert \C \rvert-1)(4\lambda_C+58)$, there exists a proper arrangement $\tau^*$ of $\Se^*$ such that the skeleton of $\tau^*$ with respect to $(\Se_1^*,\Se_2^*)$ is 2-cell and has a respectful tangle $\T^*$ of order at least $\theta^*+\phi(\rho^*)+2\rho^*$ by repeatedly applying Lemma \ref{extend a vortex} $\lvert \C \rvert$ times.
Note that $G^*$ is a subgraph of $G'$ and hence is a minor of $G$.
So Conclusion 4 holds.
Since $G^*$ is a subgraph of $G'$, Conclusion 5 holds.
For any $x \in S_x, y \in S_y$ for different members $(S_x,\Omega_x),(S_y,\Omega_y) \in \Se_2^*$, let $\Lambda^*_x,\Lambda^*_y$ be the zones $\Lambda^*$ containing $S_x$ and $S_y$, respectively, then $m_{\T^*}(x,y) \geq m_{\T'}(x,y)-\lvert \C \rvert(4\lambda_\C+58) \geq \phi(\rho^*)$.
This proves Conclusion 6.

It remains to prove that $\Se^*$ is a $\T$-central segregation of $G$.
Since $\T'$ has order at least $\theta$ and is conformal with $\T$, the order of $\T$ is at least $\theta$.
Since $\Se^*_1 \subseteq \Se_1$ and $\Se$ is $\T$-central, by Lemma \ref{make segregation central}, it is sufficient to show that there is no $(A,B) \in \T$ of order at most $2\rho^*+1$ such that $B \subseteq S$ for some $(S,\Omega) \in \Se_2^*$.
Suppose that such $(A,B)$ exists.
Since $B \subseteq S$, $V(B \cap G^*)$ consists of a subset of $V(A) \cap V(B)$ and some vertices of degree at most two in $G^*$.
Let $\beta$ be the $G^*$-minor of $G$ witnessing that $\T^*$ is conformal with $\T$.
Since $\T^*$ has order at least $\rho^*+1$, there exists a separation $(A',B') \in \T^*$ such that $\beta(E(A')) = E(A) \cap \beta(E(G^*))$, and subject to this, $\lvert V(B') \rvert$ is minimum.
Note that every vertex in $V(B')-V(A')$ has degree at most two in $G^*$.
So if $V(B')-V(A') \neq \emptyset$, then one can move a path in $B'$ with at least one internal vertex in $V(B')-V(A')$ into $A'$, contradicting the minimality of $\lvert V(B') \rvert$.
So $V(B') \subseteq V(A) \cap V(B)$ contains at most $2\rho^*+1$ vertices.
This implies that $(G^*-E(B'), B') \in \T^*$ by the second tangle axiom, contradicting the third tangle axiom.
Hence $\Se^*$ is $\T$-central.
\end{pf}

\bigskip

A segregation $\Se$ of $G$ is {\it maximal} if there exists no segregation $\Se'$ such that $\{(S,\Omega) \in \Se: \lvert \bar{\Omega} \rvert > 3\} = \{(S',\Omega') \in \Se': \lvert \overline{\Omega'} \rvert > 3\}$ and for every $(S,\Omega) \in \Se$ with $2 \leq \lvert \bar{\Omega} \rvert \leq 3$, there exists $(S',\Omega') \in \Se'$ with $\lvert \overline{\Omega'} \rvert \leq 3$ such that $S' \subseteq S$, and the containment is strict for at least one society.
Note that if a segregation $\Se$ of $G$ is maximal, then $G$ contains the skeleton of $\Se$ as a minor, and for every $(S,\Omega) \in \Se$ with $\lvert \bar{\Omega} \rvert \leq 3$ and $x \in \bar{\Omega}$, there exist $\lvert \bar{\Omega} \rvert-1$ paths in $S$ from $x$ to $\bar{\Omega}-\{x\}$ intersecting only in $\{x\}$.
Consequently, if $H$ is a triangle-free graph and the skeleton of a maximal segregation $\Se$ of $G$ admits an $H$-subdivision, then $G$ admits an $H$-subdivision.

The following theorem is proved in \cite{d} and is a stronger form of the structure theorem for excluding minors in \cite{rs XVI}.
(We remark that our definition of maximal segregations is slightly different from the one in \cite{d}.
But our alternative is required, for otherwise the skeleton mentioned in Theorem \ref{stronger excluding minor} might not be 2-cell and the metric $m_{\T'}$ might not be defined.
And the proof of Theorem \ref{stronger excluding minor} in \cite{d} works under this setting.)

\begin{theorem}[{\cite[Theorem 7]{d}}] \label{stronger excluding minor}
For every graph $L$, there exists an integer $\kappa$ such that for any nondecreasing positive function $\phi$, there exist integers $\theta, \xi, \rho$ with the following property.
Let $\T$ be a tangle of order at least $\theta$ in a graph $G$ controlling no $L$-minor of $G$.
Then there exist $Z \subseteq V(G)$ with size at most $\xi$ and a maximal $(\T-Z)$-central segregation $\Se=\Se_1 \cup \Se_2$ of $G-Z$ properly arranged by an arrangement $\alpha$ in a surface $\Sigma$ in which $L$ cannot be drawn, where every $(S,\Omega) \in \Se_1$ has the property that $\lvert \bar{\Omega} \rvert \leq 3$, and $\lvert \Se_2 \rvert \leq \kappa$ and there exists $p \leq \rho$ such that every member in $\Se_2$ is a $p$-vortex.
Furthermore, the skeleton $G'$ of $\alpha$ with respect to $(\Se_1,\Se_2)$ is $2$-cell embedded in $\Sigma$ with a respectful tangle $\T'$ of order at least $\phi(p)$ conformal with $\T-Z$, and if $x$ and $y$ are two vertices in $G'$ incident with two different members in $\Se_2$, then $m_{\T'}(x,y) \geq \phi(p)$.
\end{theorem}

Let us recall that the function mf was defined prior to Theorem \ref{main}.
A graph $H$ has a {\it nice} embedding in $\Sigma$ if $H$ can be $2$-cell embedded in $\Sigma$ and it has a set $F$ of regions such that every vertex of $H$ of degree at least $4$ is incident with exactly one region in $F$, and $\lvert F \rvert = \mf(H,\Sigma)$.

\begin{lemma}[{\cite[Lemma 12]{d}}] \label{nice H}
Let $H$ be a graph of maximum degree $d$ that can be embedded in a surface $\Sigma$.
Then there exists a triangle-free graph $H'$ of maximum degree $d$ admitting an $H$-subdivision such that $\mf(H',\Sigma) = \mf(H,\Sigma)$ and $H'$ has a nice embedding in $\Sigma$.
\end{lemma}

Recall that a vertex $v$ in a graph $G$ is $d$-free with respect to a tangle $\T$ in $G$ if there does not exist a separation $(A,B) \in \T$ of order less than $d$ such that $v \in V(A)-V(B)$.
Now, we are ready to prove Theorem \ref{stronger excluding subdivision}, which is the main theorem of this paper.
Note that Theorem \ref{main} immediately follows if we take $X=V(G)$ in Theorem \ref{stronger excluding subdivision}.

\begin{theorem} \label{stronger excluding subdivision}
Let $d,h$ be positive integers.
Then there exist $\theta, \kappa, \rho, \xi, g \geq 0$ satisfying the following properties.
If $d \geq 4$, $H$ is a graph of maximum degree at most $d$ on $h$ vertices, $G$ is a graph, and $X$ is a subset of $V(G)$ such that $G$ does not admit an $H$-subdivision whose branch vertices corresponding to vertices of degree at least four in $H$ are contained in $X$, then for every tangle $\T$ in $G$ of order at least $\theta$, there exists $Z \subseteq V(G)$ with $\lvert Z \rvert \leq \xi$ such that either
		\begin{enumerate}
			\item no vertex in $V(G-Z) \cap X$ is $d$-free with respect to $\T-Z$, or
			\item there exist a $(\T-Z)$-central segregation $\Se = \Se_1 \cup \Se_2$ of $G-Z$ with $\lvert \Se_2 \rvert \leq \kappa$, having a proper arrangement in some surface $\Sigma$ of genus at most $g$ such that every society $(S_1,\Omega_1)$ in $\Se_1$ satisfies that $\lvert \overline{\Omega_1} \rvert \leq 3$, and every society $(S_2, \Omega_2)$ in $\Se_2$ is a $\rho$-vortex, and satisfies the following property: either
			\begin{enumerate}
				\item $H$ cannot be drawn in $\Sigma$, or
				\item $H$ can be drawn in $\Sigma$ and $\mf(H,\Sigma) \geq 2$, and there exists $\Se'_2 \subseteq \Se_2$ with $\lvert \Se'_2 \rvert \leq \mf(H,\Sigma)-1$ such that every $d$-free vertex in $V(G-Z) \cap X$ with respect to $\T-Z$ is in $S-\bar{\Omega}$ for some $(S,\Omega) \in \Se'_2$.
			\end{enumerate}
		\end{enumerate}
Furthermore, if $\T$ controls a $K_{\lfloor\frac{3}{2}dh \rfloor}$-minor and $G$ does not admit an $H$-subdivision with branch vertices contained in $X$, then the first conclusion always holds even when $d \leq 3$.
\end{theorem}

\begin{pf}
Note that there are only finitely many graphs of maximum degree at most $d$ on $h$ vertices, and there are only finitely many surfaces in which $H$ can be drawn but $K_{\lfloor\frac{3}{2}dh \rfloor}$ cannot.
So there exists $h^*$ such that for every graph $H$ on $h$ vertices of maximum degree at most $d$ and surface in which $H$ can be drawn but $K_{\lfloor\frac{3}{2}dh \rfloor}$ cannot, the graph $H'$ mentioned in Lemma \ref{nice H} can be chosen such that it has at most $h^*$ vertices.

We define the following.
	\begin{itemize}
		\item Let $\kappa_{\ref{stronger excluding minor}}$ be the number $\kappa$ mentioned in Theorem \ref{stronger excluding minor} by taking $L=K_{\lfloor\frac{3}{2}dh \rfloor}$.
		\item Let $\theta_{\ref{free in vortex}}, \beta_{\ref{free in vortex}}, f_{\ref{free in vortex}}$ be the functions $\theta_0, \beta, f$ mentioned in Lemma \ref{free in vortex}, respectively.
		\item Let $\phi'$ be the maximum of $\phi_{\ref{free should far apart}}(d,h^*,\Sigma)$ among all surfaces $\Sigma$ in which $K_{\lfloor\frac{3}{2}dh \rfloor}$ cannot be drawn, where $\phi_{\ref{free should far apart}}$ is the number $\phi$ mentioned in Lemma \ref{free should far apart}.
		\item Let $\theta_{\ref{linkage on surface}}$ be the maximum of $\theta$ mentioned in Theorem \ref{linkage on surface} by taking all surfaces in which $K_{\lfloor\frac{3}{2}dh \rfloor}$ cannot be drawn, $t=(d+2)h^*$ and $z=dh^*$.
		\item Let $\phi^*(x) = \theta_{\ref{free in vortex}}(d,h^*,x, 2\kappa_{\ref{stronger excluding minor}}+h^*,(dh^*+h^*+1)(\theta_{\ref{linkage on surface}}+1))+2f_{\ref{free in vortex}}(d,h^*,x,2\kappa_{\ref{stronger excluding minor}}+h^*) + (2\kappa_{\ref{stronger excluding minor}}+h^*)(6\beta_{\ref{free in vortex}}(d,h^*,x)+4f_{\ref{free in vortex}}(d,h^*,x,2\kappa_{\ref{stronger excluding minor}}+h^*)+2) +2(dh^*+h^*+1)(\theta_{\ref{linkage on surface}}+6)$.
		\item Let $\theta_{\ref{make vortex away}}'(x)$ be the function $\theta$ obtained by applying Lemma \ref{make vortex away} by taking $\phi = \phi^*$, $\rho=x$, $\lambda=d+\phi'+11$, $\kappa=\kappa_{\ref{stronger excluding minor}}$, $k=h^*+\kappa_{\ref{stronger excluding minor}}$, $\theta^*=\theta_{\ref{free in vortex}}(d,h^*,x, 2\kappa_{\ref{stronger excluding minor}}+h^*,(dh^*+h^*+1)(\theta_{\ref{linkage on surface}}+1))$ and $d=d$. 
		\item Let $\theta_{\ref{free should far apart}}$ be the maximum of $\theta(d,h^*,\Sigma)$ mentioned in Lemma \ref{free should far apart} among all surfaces $\Sigma$ in which $K_{\lfloor\frac{3}{2}dh \rfloor}$ cannot be drawn.
		\item Let $\theta_{\ref{stronger excluding minor}}$, $\xi_{\ref{stronger excluding minor}}$, $\rho_{\ref{stronger excluding minor}}$ be the numbers $\theta, \xi, \rho$ mentioned in Theorem \ref{stronger excluding minor}, respectively, by taking $L=K_{\lfloor\frac{3}{2}dh \rfloor}$ and further taking $\phi(x)=\theta_{\ref{make vortex away}}'(x)+\theta_{\ref{free should far apart}}$.
		\item Let $\theta_{\ref{make vortex away}}$ and $\rho_{\ref{make vortex away}}$ be the numbers $\theta$ and $\rho^*$ obtained by applying Lemma \ref{make vortex away} by taking $\phi= \phi^*$, $\rho=\rho_{\ref{stronger excluding minor}}$, $\lambda=d+\phi'+11$, $\kappa=\kappa_{\ref{stronger excluding minor}}$, $k=h^*+\kappa_{\ref{stronger excluding minor}}$, $\theta^*=\theta_{\ref{free in vortex}}(d,h^*,\rho_{\ref{stronger excluding minor}}, 2\kappa_{\ref{stronger excluding minor}}+h^*,(dh^*+h^*+1)(\theta_{\ref{linkage on surface}}+1))$ and $d=d$.
		\item Let $\theta_{\ref{spider tangle}}=(hd)^{d+1}+d$. 
	\end{itemize}

Now we are ready to define the numbers for the conclusion of this theorem.
	\begin{itemize}
		\item Let $\xi = \max\{\xi_{\ref{stronger excluding minor}}+(2\kappa_{\ref{stronger excluding minor}}+h^*)\beta_{\ref{free in vortex}}(d,h^*,\rho_{\ref{make vortex away}}) , (hd+1)^{d+1}\}$.
		\item Let $\theta=\theta_{\ref{make vortex away}}+ \theta_{\ref{free should far apart}} + \theta_{\ref{stronger excluding minor}} + \xi+d+\rho_{\ref{stronger excluding minor}}+\rho_{\ref{make vortex away}}$. 
		\item Let $\kappa=\kappa_{\ref{stronger excluding minor}}+h^*$.
		\item Let $\rho=\rho_{\ref{stronger excluding minor}}+\rho_{\ref{make vortex away}}$.
		\item Let $g$ be the maximum genus of a surface in which $K_{\lceil \frac{3}{2}dh \rceil}$ cannot be drawn.
	\end{itemize}

Let $\T$ be a tangle of order at least $\theta$ in $G$, and assume that $G$ has no $H$-subdivision with branch vertices contained in $X$.
We may assume that $X$ contains at least $h$ vertices of degree at least $d$ in $G$, otherwise the first statement holds by letting $Z$ be the set of vertices in $X$ of degree at least $d$.
We first assume that $\T$ controls a $K_{\lfloor\frac{3}{2}dh \rfloor}$-minor.
By Lemma \ref{vertex spider subdivision} and Theorem \ref{spider tangle}, since $G$ does not admit an $H$-subdivision with branch vertices contained in $X$, there exists a set of vertices $Z$ of $G$ with $\lvert Z \rvert \leq \xi$ such that for every vertex $v \in V(G-Z) \cap X$ of degree at least $d$ in $G$, there exists a separation $(A_v,B_v) \in \T-Z$ of $G-Z$ of order at most $d-1$ such that $v \in V(A_v)-V(B_v)$.
Therefore, the first statement holds.

So we may assume that $\T$ does not control a $K_{\lfloor\frac{3}{2}dh \rfloor}$-minor.
To prove this theorem, we may assume that $G$ does not contain an $H$-subdivision whose branch vertices corresponding to vertices of degree at least four in $H$ are contained in $X$.

By Theorem \ref{stronger excluding minor}, there exist a surface $\Sigma$ in which $K_{\lfloor\frac{3}{2}dh \rfloor}$ cannot be drawn, $Z \subseteq V(G)$ with $\lvert Z \rvert \leq \xi_{\ref{stronger excluding minor}}$, a number $p$ with $p \leq \rho_{\ref{stronger excluding minor}}$, and a maximal $(\T-Z)$-central segregation $\Se = \Se_1 \cup \Se_2$ of $G-Z$ with $\lvert \Se_2 \rvert \leq \kappa_{\ref{stronger excluding minor}}$, having a proper arrangement $\tau$ in $\Sigma$ such that every society $(S,\Omega)$ in $\Se_1$ satisfies that $\lvert \bar{\Omega} \rvert \leq 3$, and every society in $\Se_2$ is a $p$-vortex, and the skeleton $G'$ of $\Se$ is $2$-cell embedded in $\Sigma$ and has a respectful tangle $\T'$ of order at least $\theta_{\ref{make vortex away}}'(p)+\theta_{\ref{free should far apart}}$ conformal with $\T-Z$, and if $x,y$ are two vertices in $G'$ incident with two different members in $\Se_2$, then $m_{\T'}(x,y) \geq \theta_{\ref{make vortex away}}'(p)$.
If $H$ cannot be drawn in $\Sigma$, then Statement 2(a) holds, so we may assume that $H$ can be drawn in $\Sigma$.

In addition, we may assume that $V(G-Z) \cap X$ contains $d$-free vertices with respect to $\T-Z$, for otherwise Statement 1 holds.
Note that every vertex in $\bigcup_{(S,\Omega) \in \Se_1}V(S)- V(G')$ is not $d$-free with respect to $\T-Z$ since $d \geq 4$.
If $v$ is in $V((G-Z) \cap G')$ but is not $d$-free with respect to $\T'$, then there exists a separation $(A',B') \in \T'$ of order less than $d$ such that $v \in V(A')-V(B')$.
We choose $(A',B')$ such that $A'$ is as small as possible.
Note $m_{\T'}(v,x) <d$ for every $x \in V(A')$ by Theorem \ref{A distance}.
Suppose that there is no vertex $x \in \overline{\Omega}$ with $(S,\Omega) \in \Se_2$ and $m_{\T'}(v,x) < d$.
Then there exists $(A,B) \in \T-Z$ of order less than $d$ such that $V(A) = \bigcup_{(S,\Omega) \in \Se, \overline{\Omega} \subseteq V(A')} V(S)$ and $V(A) \cap V(B) = V(A') \cap V(B')$.
So $v$ is not $d$-free with respect to $\T-Z$.
Therefore, if $v$ is a vertex in $(G-Z) \cap G'$ that is $d$-free with respect to $\T-Z$ but not $d$-free with respect to $\T'$, then $m_{\T'}(v,x) < d$ for some $x \in V(S)$ with $(S,\Omega) \in \Se_2$.
By Theorem \ref{big zone contains ball} and Lemma \ref{disjoint boundry of zone}, for every $(S,\Omega) \in \Se_2$, there exists a $(d+11)$-zone $\Lambda_S$ with respect to $\T'$ around a vertex in $\bar{\Omega}$ containing every atom $y$ with $m_{\T'}(x,y) \leq d+1$ for one such $x$.
Thus every vertex of $(G-Z) \cap G'$ that is $d$-free with respect to $\T-Z$ but not $d$-free with respect to $\T'$ is in $\bigcup_{(S,\Omega) \in \Se_2} \Lambda_S$.

Let $H'$ be a graph that has a nice embedding mentioned in Lemma \ref{nice H} such that $\lvert V(H') \rvert \leq h^*$.
By Lemma \ref{free should far apart}, $V(G') \cap X$ does not contain $\lvert V(H') \rvert$ $d$-free vertices with respect to $\T'$ such that every pair of them has distance at least $\phi'$ under the metric $m_{\T'}$, for otherwise $G$ contains an $H$-subdivision with branch vertices contained in $X$.
So by Theorem \ref{big zone contains ball} and Lemma \ref{disjoint boundry of zone}, there exist integer $k$ with $0 \leq k \leq h^*$, $d$-free vertices $v_1, v_2,...,v_k$ of $G'$ with respect to $\T'$, and $(\phi'+10)$-zones $\Lambda_1, \Lambda_2, ..., \Lambda_k$ around $v_1,v_2,...,v_k$, respectively, such that every $d$-free vertex in $V(G') \cap X$ with respect to $\T'$ is in $\bigcup_{i=1}^k \Lambda_i$.
Then every $d$-free vertex in $V(G-Z) \cap X$ with respect to $\T-Z$ is a vertex of $G' \cup \bigcup_{(S,\Omega)\in\Se_2}S$, and it is in $\bigcup_{i=1}^k \Lambda_i \cup \bigcup_{(S,\Omega) \in \Se_2} \Lambda_S$.

Then let $\Se^*=\Se^*_1 \cup \Se_2^*$, $\rho',\T^*$ and $G^*$ be the $\Se^*$, $\rho',\T^*$ and $G^*$, respectively, mentioned in the conclusion of Lemma \ref{make vortex away} by taking $\phi=\phi^*$, $\rho=p$, $\lambda=d+\phi'+11$, $\kappa=\kappa_{\ref{stronger excluding minor}}$, $k=h^*+\kappa_{\ref{stronger excluding minor}}$, $\theta^*=\theta_{\ref{free in vortex}}(d,h^*,\rho_{\ref{make vortex away}}, 2\kappa_{\ref{stronger excluding minor}}+h^*,(dh^*+h^*+1)(\theta_{\ref{linkage on surface}}+1))$ and $d=d$, and further taking $G=G-Z$, $\T=\T-Z$, $\Se=\Se$, $\tau=\tau$, $\Sigma=\Sigma$, and $G'$ to be the skeleton of $\Se$.
Note that $\lvert \Se_2^* \rvert \leq 2\kappa_{\ref{stronger excluding minor}}+h^*$, and the order of $\T^*$ is at least $\theta_{\ref{free in vortex}}(d,h^*,\rho_{\ref{make vortex away}}, 2\kappa_{\ref{stronger excluding minor}}+h^*,(dh^*+h^*+1)(\theta_{\ref{linkage on surface}}+1))+\phi^*(\rho_{\ref{make vortex away}})+2\rho_{\ref{make vortex away}} \geq \theta_{\ref{free in vortex}}(d,h^*,\rho_{\ref{make vortex away}}, 2\kappa_{\ref{stronger excluding minor}}+h^*,(dh^*+h^*+1)(\theta_{\ref{linkage on surface}}+1))$.

Let $\kappa'$ be the number of members of $\Se_2^*$ containing $d$-free vertices with respect to $\T-Z$ belonging to $X$.
Note that $\kappa' \leq \lvert \Se_2^* \rvert \leq 2\kappa_{\ref{stronger excluding minor}}+h^*$.
Let $Z_1,Z_2,...,Z_{\kappa'}, \Lambda_1,\Lambda_2,...,\Lambda_{\kappa'}$ be the sets obtained by applying Lemma \ref{free in vortex} by taking $\kappa=\kappa'$, $h_i=h^*$ for every $i$, $\rho = \rho'$, $\theta'' = (dh^*+h^*+1)(\theta_{\ref{linkage on surface}}+1)$, $G=G-Z$, $G'=G^*$ and $(S_1,\Omega_1),(S_2,\Omega_2),...,(S_{\kappa'},\Omega_{\kappa'})$ as the vortices in $\Se_2^*$ containing $d$-free vertices with respect to $\T-Z$ belonging to $X$.
Define ${\Se_2^*}' \subseteq \Se_2^*$ to consist of the members in which $\Lambda_i \neq \emptyset$.
We define $Z'$ to be $Z \cup \bigcup_{1 \leq i \leq \kappa'} Z_i$.
Note that $\lvert Z' \rvert \leq \xi$.

We may assume that there exist $d$-free vertices of $V(G-Z') \cap X$ with respect to $\T-Z'$; otherwise Statement 1 holds.
So $\lvert {\Se_2^*}' \rvert \geq 1$.
In addition, we may assume that $\mf(H,\Sigma) \geq 1$ since otherwise $H$ contains no vertex of degree at least four, and hence $G$ has no $H$ minor and Statement 2(a) holds by Theorem \ref{stronger excluding minor}.
If $\mf(H,\Sigma) \geq 2$ and $\lvert {\Se_2^*}' \rvert \leq \mf(H,\Sigma)-1$, then Statement 2(b) holds.
So we may assume that $\lvert {\Se_2^*}' \rvert \geq \mf(H,\Sigma) \geq 1$ and that $\Lambda_i \neq \emptyset$ for $i=1,2,...,\lvert {\Se_2^*}' \rvert$.

Let $G''$ be the drawing and $\T''$ the tangle in $G''$ conformal with $\T^*$ mentioned in the conclusion of Lemma \ref{free in vortex}.
For $1 \leq i \leq \lvert {\Se_2^*}' \rvert$ and $1 \leq j \leq h^*$, let $Y_i$ and $A_{i,j}$ be the cycles and sets mentioned in Conclusion 2 of Lemma \ref{free in vortex}, respectively.
So for every $1 \leq i < i' \leq \lvert {\Se_2^*}' \rvert$, $j,j' \in \{1,2,...,h^*\}$, $x \in A_{i,j},y \in A_{i',j'}$, we have that $m_{\T''}(x,y) \geq m_{\T^*}(x,y)- (2\kappa_{\ref{stronger excluding minor}}+h^*)(4f_{\ref{free in vortex}}(d,h^*,\rho',2\kappa_{\ref{stronger excluding minor}}+h^*)+2)$.

Note that for such $x,y$, we know $m_{\T^*}(x,y) \geq \phi^*(\rho')-2f_{\ref{free in vortex}}(d,h^*,\rho', 2\kappa_{\ref{stronger excluding minor}}+h^*)$, since each $x,y$ is within distance at most $f_{\ref{free in vortex}}(d,h^*,\rho', 2\kappa_{\ref{stronger excluding minor}}+h^*)$ (with respect to $m_{\T^*}$) away from a member of ${\Se_2^*}'$, and the distance (with respect to $m_{\T^*}$) between those two members is at least $\phi^*(\rho')$.
Therefore, for every $1 \leq i < i' \leq \lvert {\Se_2^*}' \rvert$, $j,j' \in \{1,2,...,h^*\}$, $x \in A_{i,j},y \in A_{i',j'}$, we have that $m_{\T''}(x,y) \geq \phi^*(\rho')-2f_{\ref{free in vortex}}(d,h^*,\rho', 2\kappa_{\ref{stronger excluding minor}}+h^*)-(2\kappa_{\ref{stronger excluding minor}}+h^*)(4f_{\ref{free in vortex}}(d,h^*,\rho',2\kappa_{\ref{stronger excluding minor}}+h^*)+2) \geq 2(dh^*+h^*+1)(\theta_{\ref{linkage on surface}}+6)$. 

Let $x \in A_{1,1}$.
By Lemma \ref{something is far}, there exists an edge $e^*$ of $G''$ with $m_{\T''}(e^*,x) \geq (dh^*+h^*+1)(\theta_{\ref{linkage on surface}}+6)$.
As in the proof of Theorem 4.3 in \cite{rs XII}, there exist non-loop edges $e_1,e_2,...,e_{dh^*+h^*}$ of $G''$ such that $(\theta_{\ref{linkage on surface}}+6)i \leq m_{\T''}(x,e_i) \leq (\theta_{\ref{linkage on surface}}+6)i+3$ for $1 \leq i \leq dh^*+h^*$, and the set of the ends of each $e_i$ is free for each $1 \leq i \leq dh^*+h^*$.
Therefore, $m_{\T''}(e_{i},e_j) \geq \theta_{\ref{linkage on surface}}+3$ for every $1 \leq i<j \leq dh^*+h^*$.
Note that $m_{\T''}(x,y) \leq 2$ for $y \in \bigcup_{j=1}^{h^*}A_{1,j}$ and $m_{\T''}(x,y) \geq 2(dh^*+h^*+1)(\theta_{\ref{linkage on surface}}+6)$ for $y \in \bigcup_{i=2}^{\lvert {\Se_2^*}' \rvert} \bigcup_{j=1}^{h^*}A_{i,j}$.
Hence, $m_{\T''}(y,e_\ell) \geq \theta_{\ref{linkage on surface}}+1$ for every $y \in \bigcup_{i=1}^{\lvert {\Se_2^*}' \rvert} \bigcup_{j=1}^{h^*} A_{i,j}$ and $1 \leq \ell \leq dh^*+h^*$.
For $1 \leq i \leq \lvert {\Se_2^*}' \rvert$, define $\Delta_i$ to be a disk in $\Sigma$ contained in $\overline{\Lambda_i}$ such that $\Delta_i \cap G'' = \bigcup_{j=1}^{h^*}A_{i,j}$.
For $1 \leq i \leq dh^*+h^*$, define $\Delta_{\lvert {\Se_2^*}' \rvert+i}$ to be a disk in $\Sigma$ such that $\Delta_i \cap G''$ is the set of the ends of $e_i$.
Let $W=\bigcup_{i=1}^{\mf(H',\Sigma)}\bigcup_{j=1}^{h^*} A_{i,j} \cup \bigcup_{i=1}^{\lvert E(H') \rvert+\lvert V(H') \rvert}\{a_i,b_i\}$, where $a_i,b_i$ are the ends of $e_i$.
By Theorem \ref{linkage on surface}, for any positive integer $p$ and any partition $(W_1,...,W_p)$ of $W$ satisfying the topological feasibility condition, there exist pairwise disjoint connected subdrawings $\Gamma_1,...,\Gamma_p$ of $G''$ with $V(\Gamma_j) \cap W=W_j$ for every $1 \leq j \leq p$.

Since $H'$ has a nice embedding in $\Sigma$, we can embed $H'$ into $\Sigma$ such that the vertices of degree at least four of $H'$ are incident with $\mf(H',\Sigma)$ regions.
Then there exists a partition $(W_z: z \in V(H') \cup E(H'))$ (with possible empty part) of a subset of $W$ such that the following hold.
	\begin{itemize}
		\item Fix an nice embedding of $H'$ in $\Sigma$, a set of $\mf(H',\Sigma)$ faces incident with all vertices of degree at least four in $H'$, an injection $\iota_1$ that maps each of those $\mf(H',\Sigma)$ faces to a member of ${\Se_2^*}'$, a set of $h^*$ disjoint $d$-spiders mentioned in Conclusion 2(d) in Lemma \ref{free in vortex} for each member of ${\Se_2^*}'$, and an injection $\iota_2$ that maps each vertex $u$ of $H'$ of degree at least four in $H'$ to a $d$-spider in the set for $\iota_1(f)$, where $f$ is the face among those $\mf(H',\Sigma)$ faces incident with $u$.
		\item There exist an injection $\iota_3$ from $E(H')$ to $\{e_t: 1 \leq t \leq dh^*\}$, and an injection $\iota_4$ that maps each vertex-edge incidence pair $(v,e)$ of $H'$ to an end of $\iota_3(e)$, where every loop contributes two vertex-edge incidence pairs.
		\item For each vertex $v$ of $H'$ of degree at most three, $W_v=\{\iota_4((v,e)): e$ is an edge of $H'$ incident with $v\}$.
		\item For each edge $e$ of $H'$ incident with both ends, say $u,v$, of degree at least four, where $u,v$ are not necessarily distinct, $W_e$ is a set of size two consisting of one leaf of $\iota_2(u)$ and one leaf of $\iota_2(v)$.
		\item For each edge $e$ of $H'$ incident with one vertex, say $u$, of degree at most three and one vertex, say $v$, of degree at least four, $W_e$ consists of $\iota_3((u,e))$ and one leaf of $\iota_2(v)$.
		\item $W_z=\emptyset$ if $z$ is a vertex of $H'$ of degree at least four or $z$ is an edge of $H'$ not incident with any vertex of degree at least four.
		\item $(W_z: z \in V(H') \cup E(H'))$ satisfies the topological feasibility condition.
	\end{itemize}
Since the partition $(W_z: z \in V(H') \cup E(H'))$ satisfies the topological feasibility condition, by Theorem \ref{linkage on surface}, there exist pairwise disjoint connected subdrawings $(\Gamma_z: z \in V(H') \cup E(H'))$ of $G''$ with $V(\Gamma_z) \cap W=W_z$ for every $z \in V(H') \cup E(H')$.
Observe that the union of the image of $\iota_2$ and $\bigcup_{z \in V(H') \cup E(H')}\Gamma_z \cup \bigcup_{i=1}^{dh^*}e_i$ contains an $H'$-subdivision $(\pi_V,\pi_E)$ such that its every branch vertex corresponding to a vertex of degree at least four in $H'$ is contained in $X$.

Finally, we shall obtain a contradiction by showing that $G$ admits an $H$-subdivision whose branch vertices corresponding to vertices of degree at least four in $H$ are contained in $X$.
Recall that $\Se^*$ is maximal, so for every $(S,\Omega) \in \Se^*_1$ and for every $a \in \bar{\Omega}$, there exist $\lvert \bar{\Omega} \rvert-1$ paths in $S$ from $a$ to $\bar{\Omega}-\{a\}$ intersecting in $a$ and otherwise disjoint.
Since $H'$ is triangle-free, one edge of the triangle in $G''$ formed by $\bar{\Omega}$ is not contained in the image of $\pi_E$ for any $(S,\Omega) \in \Se^*_1$ with $\lvert \bar{\Omega} \rvert=3$.
Therefore, $G$ admits an $H$-subdivision whose branch vertices corresponding to vertices of degree at least four in $H$ are contained in $X$.
This completes the proof.
\end{pf}

\section{Remarks about optimality}
\label{sec:bestmaxdeg}

In this section we shall prove that the order of the separations mentioned in the first conclusion of Theorem \ref{main} cannot be improved, and the value $D$ mentioned in Theorem \ref{gm structure} and \cite[Theorem 3]{d} cannot be replaced by any number smaller than the number $d-1$ mentioned in the first conclusion of Theorem \ref{main}.

Let $\theta$ be a positive integer.
A {\it $2\theta \times 2\theta$-wall} is the graph $W$ with $V(W)=\{(i,j): 1 \leq i \leq 2\theta, 1 \leq j \leq 2\theta\}$ and $E(W)=\{(i,j)(i,j+1): 1 \leq i \leq 2\theta, 1 \leq j \leq 2\theta-1\} \cup \{(2k-1,j)(2k,j): 1 \leq k \leq \theta, 1 \leq j \leq 2\theta$, $j$ is odd$\} \cup \{(2k,j)(2k+1,j): 1 \leq j \leq \theta-1, 1 \leq j \leq 2\theta$, $j$ is even$\}$.
For each $1 \leq i \leq 2\theta$, the {\it $i$-th row} is the subgraph induced by the vertices $\{(i,j): 1 \leq j \leq 2\theta\}$.

We fix $d$ to be a positive integer with $d \geq 4$ in the rest of this section.
Let $r$ be a positive integer.
Let $Y$ be the collection of all $(d-1)$-element subsets of the edge-set of the $2r \times 2r$-wall.
Define $G_r$ to be a graph obtained from the $2r \times 2r$-wall by the following procedure.
	\begin{itemize}
		\item Subdividing each edge $\lvert Y \rvert$ times.
		\item For each member $y$ of $Y$, add a vertex $v_y$.
			Let $V_Y=\{v_y: y \in Y\}$.
		\item For each $y \in Y$ and $e \in y$, add an edge between $v_y$ and a vertex obtained by subdividing $e$ in a way that every vertex not in $V_Y$ has degree at most three in $G_r$.
	\end{itemize}
Note that $G_r$ has maximum degree $d-1$, so $G_r$ does not contain any graph with maximum degree at least $d$ as a subdivision.

\begin{theorem}
For all positive integers $d,\theta,\kappa,\rho,\xi,g$, there exist a positive integer $r$ and a tangle $\T$ in $G_r$ of order at least $\theta$ such that the following hold.
	\begin{enumerate}
		\item For every $Z \subseteq V(G_r)$ with $\lvert Z \rvert \leq \xi$, there exists a vertex $v \in V_Y-Z$ such that there exists no separation $(A,B) \in \T-Z$ of order at most $d-2$ with $v \in V(A)-V(B)$.
		\item For every $Z \subseteq V(G_r)$ with $\lvert Z \rvert \leq \xi$, there exists no $(\T-Z)$-central segregation $\Se=\Se_1 \cup \Se_2$ of $G_r-Z$ with $\lvert \Se_2 \rvert \leq \kappa$ such that $\Se$ can be properly arranged in some surface $\Sigma$ with genus at most $g$, every $(S_1,\Omega_1) \in \Se_1$ satisfies $\lvert \overline{\Omega_1} \rvert \leq 3$, and every $(S_2,\Omega_2) \in \Se_2$ is a $\rho$-vortex.
		\item $G_r$ cannot be constructed by clique-sums, starting with graphs that are an $\leq \xi$-extension of either a graph of maximum degree at most $d-2$, or an out-growth by $\leq \kappa$ $\rho$-rings of a graph that can be drawn in a surface of genus at most $g$.
	\end{enumerate}
\end{theorem}

\begin{pf}
We may assume $d \geq 2$; otherwise the first conclusion holds.
Let $r$ be a fixed positive integer.
We shall prove that every sufficiently large $r$ satisfies the conclusion of this theorem.

Define $\T$ to be the collection of all separations $(A,B)$ of $G_r$ of order less than $\lfloor \frac{2r}{3} \rfloor$ such that $B$ contains all vertices in some row of the $2r \times 2r$-wall.
It is straightforward to verify that $\T$ is a tangle of order $\lfloor \frac{2r}{3} \rfloor$.
We choose $r$ to be a number with $\lfloor \frac{2r}{3} \rfloor >\theta$, then $\T$ is a tangle in $G_r$ of order at least $\theta$.

There exists a positive integer $k$ such that if $\T$ controls a $K_k$-minor, then there exist no $Z \subseteq V(G_r)$ with $\lvert Z \rvert \leq \xi$ and a $(\T-Z)$-central segregation $\Se=\Se_1 \cup \Se_2$ of $G_r-Z$ with $\lvert \Se_2 \rvert \leq \kappa$ such that $\Se$ can be properly arranged in a surface of genus at most $g$, every $(S_1,\Omega_1) \in \Se_1$ satisfies $\lvert \Omega_1 \rvert \leq 3$, and every $(S_2,\Omega_2) \in \Se_2$ is a $\rho$-vortex.
It is not hard to show that there exists an integer $r_k$ such that if $r \geq r_k$, then $G_r$ contains a $K_k$-minor such that each branch set intersects every row of the $2r \times 2r$-wall in $G_r$.
That is, if $r \geq r_k$, then $\T$ controls a $K_k$-minor.
We choose $r \geq r_k$, so the second conclusion of this theorem holds.

Choose $r \geq d^{2\xi+1}2^{2d^2}$.
Now we prove the first conclusion of this theorem.
Let $Z \subseteq V(G)$ with $\lvert Z \rvert \leq \xi$.
Then $G_r-Z$ contains a subgraph $R$ which is a subdivision of $G_{2^{2d^2}d}$ and only the edges contained in the subdivided wall are subdivided.
Let $v$ be a vertex of $V(R) \cap V_Y$ such that there exists a row of the subdivided wall in $G_{2^{2d^2}d}$ such that all neighbors of $v$ are vertices obtained by subdividing edges of that row.
Then for any row of the subdivided wall in $G_{2^{2d^2}d}$, there exist $d-1$ disjoint paths in the subdivided wall from the neighbors of $v$ to that row.
So if $(A,B) \in \T-Z$ with $v \in V(A)-V(B)$, then the order of $(A,B)$ is at least $d-1$.
This proves the first conclusion of this theorem.

Finally, we prove the third conclusion of this theorem.
Suppose that $G_r$ can be constructed by clique-sums starting from graphs mentioned in the third conclusion.
Then $G_r$ has a tree-decomposition $(T,\X)$ such that each torso is a graph involved in the clique-sum.
For each node $t$ of the tree, let $X_t$ be the bag at $t$.
Since the clique size of each graph involved in the clique-sum is bounded, the adhesion of $(T,\X)$ is bounded.
We choose $r$ to be larger than the adhesion of $(T,\X)$.
For each edge $tt'$ of the tree $T$, let $T_t$ and $T_{t'}$ be the components of $T-tt'$ containing $t$ and $t'$, respectively.
Then there exists a separation $(A_{tt'},B_{tt'})$ of $G_r$ such that $V(A_{tt'})=\bigcup_{t'' \in V(T_t)}X_{t''}$ and $V(B_{tt'})=\bigcup_{t'' \in V(T_{t'})}X_{t''}$.
Since the order of $\T$ is larger than the adhesion of $(T,\X)$, either $(A_{tt'},B_{tt'}) \in \T$ or $(B_{tt'},A_{tt'}) \in \T$.
We orient the edge $tt'$ from $t$ to $t'$ if $(A_{tt'},B_{tt'}) \in \T$.
Then there exists a node $t^*$ of $T$ with out-degree zero.
Let $L$ be the torso at $t^*$.
Since $\T$ controls a $K_k$-minor, $L$ is a $\leq \xi$-extension of a graph of maximum degree at most $d-2$.
Let $Z \subseteq V(L)$ with $\lvert Z \rvert \leq \xi$ such that $L-Z$ has maximum degree at most $d-2$.

Note that $G_r-Z$ contains a subgraph $R$ which is a subdivision of $G_{2^{2d^2}d}$ and only the edges contained in the subdivided wall are subdivided.
And there exists a vertex $v$ of $V(R) \cap V_Y$ such that there exists a row of the subdivided wall in $G_{2^{2d^2}d}$ such that all neighbors of $v$ are vertices obtained by subdividing edges of that row.

Since $L-Z$ has maximum degree at most $d-2$, $\lvert X_{t^*} \cap X_t-Z \rvert \leq d-2$ for every neighbor $t$ of $t^*$.
Since there exists no separation $(A,B) \in \T-Z$ of order at most $d-2$ with $v \in V(A)-V(B)$, $v \in X_{t^*}$.
But since $v$ has degree $d-1$ in $G_r-Z$, $v$ belongs to $X_{t^*} \cap X_c$ for some neighbor $c$ of $t^*$.
Let $c_1,c_2,...,c_\ell$ be the neighbors of $t^*$ such that $v$ is adjacent in $G_r$ to some vertex in $(\bigcup_{i=1}^\ell\bigcup_{t \in T_{c_i}}X_t)-X_{t^*}$, where $T_{c_i}$ is the component of $T-t^*c_i$ containing $c_i$.
Since $v$ has degree $d-1$ in $G_r$, $\ell \leq d-1$.
Let $(A',B')$ be a separation of $G_r$ such that $V(A')=Z \cup \bigcup_{i=1}^\ell\bigcup_{t \in V(T_{c_i})}X_t$ and $V(B')=Z \cup \bigcup_{t \in V(T)-\bigcup_{i=1}^{\ell}V(T_{c_i})}X_t$.
Note that the order of $(A',B')$ is at most $\xi+d\ell \leq \xi+d^2<\lfloor \frac{2r}{3} \rfloor$.
So $(A',B') \in \T$.
Let $(A^*,B^*) \in \T-Z$ be a separation of $G_r-Z$ with $V(A^*)=V(A')-Z$ and $V(B^*)=V(B')-Z$.
Since the order of $(A^*,B^*)$ is at most $d^2$, $R \cap B^*$ contains a subgraph $L'$ isomorphic to a subdivision of $G_{d}$, where only the edges contained in the subdivided wall are subdivided.
So there exists a row of the subdivided wall in $L'$ contained in $B^*$.
Note that there are $d-1$ disjoint paths $P_1,P_2,...,P_{d-1}$ in $G_r-Z$ from the neighbors of $v$ to that row, where each neighbor of $v$ in $G_r-Z$ is contained in some $P_i$.
Let $p$ be the number of neighbors of $v$ in $G_r-Z$ contained in $V(A^*)-V(B^*)$, and we may assume that $P_1,P_2,...,P_p$ contain such neighbors.
So for each $i$ with $1 \leq i \leq p$, $P_i$ intersects $V(A^* \cap B^*)$.
Hence $p \leq \lvert V(A^* \cap B^*)-\{v\} \rvert$.
By the definition of a torso, $v$ is adjacent in $L$ to every vertex in $V(A^* \cap B^*)-\{v\}$.
Therefore, the degree of $v$ in $L$ is at least the degree of $v$ in $G_r-Z$, which is $d-1$, a contradiction.
This proves the theorem.
\end{pf}

\bigskip

\noindent{\bf Acknowledgement:} Proofs of several results in this paper use the Graph Minors theory developed by Robertson and Seymour and are inspired by ideas in the proof of Dvo\v{r}\'ak's theorem \cite{d}.
The authors would like to acknowledge Paul Wollan for inspiring conversations.
The authors would also like to thank the anonymous referee for his or her careful reading and constructive suggestions.
The paper is based on part of the PhD dissertation~\cite{LiuPhD} of the first author.

\end{document}